\documentclass[notitlepage,twoside,a4paper]{amsart}
\usepackage{mathrsfs}
\usepackage{amsmath,amssymb,enumerate}
\usepackage{epsfig,fancyhdr,color}
\usepackage{epstopdf}
\usepackage{amssymb}
\usepackage{amsmath,amsthm}
\usepackage{latexsym}
\usepackage{amscd}
\usepackage{psfrag}
\usepackage{graphicx}
\usepackage{epsf}
\usepackage[latin1]{inputenc}
\usepackage[all]{xy}
\usepackage{tikz}
\usepackage{prettyref}
\usepackage{subfigure}
\usepackage{float}
\usepackage{color}
\usepackage{array}
\usepackage{hyperref}
\usepackage{cite}
\usepackage[T1]{fontenc}

\usepackage[totalwidth=16cm,totalheight=24cm]{geometry}

\newcommand{\III}{I\hspace{-0.1cm}I\hspace{-0.1cm}I}

\newtheorem{theorem}{\rm\bf Theorem}[section]
\newtheorem*{theorem*}{\rm\bf Theorem A}
\newtheorem*{theorem*'}{\rm\bf Theorem A*}
\newtheorem*{theorem**}{\rm\bf Theorem B}
\newtheorem*{theorem***}{\rm\bf Theorem C}
\newtheorem*{theorem****}{\rm\bf Theorem D}
\newtheorem*{question*}{\rm\bf Question}
\newtheorem*{remark*}{\rm\bf Remark}

\newtheorem{proposition}[theorem]{\rm\bf Proposition}

\newtheorem{definition}[theorem]{\rm\bf Definition}
\newtheorem{remark}[theorem]{\rm\bf Remark}

\newtheorem{question}[theorem]{\rm\bf Question}

\newtheorem{lemma}[theorem]{\rm\bf Lemma}

\newtheoremstyle{named}{}{}{\itshape}{}{\bfseries}{.}{.5em}{#1 \thmnote{#3}}
\theoremstyle{named}

\newcommand{\AdS}{\mathbb{ADS}}

\newcommand{\PSL}{\rm{PSL}}
\newcommand{\CH}{\rm{CH}}

\newcounter{notes}%

\def\interieur#1{\mathord{\mathop{\kern 0pt #1}\limits^\circ}}

\title[The Prescribed Metric on Convex Subsets of Anti-de Sitter Space with Quasi-Circle Ideal Boundaries]{The Prescribed Metric on Convex Subsets of Anti-de Sitter Space with Quasi-Circle Ideal Boundaries}

\date{v0, \today}
\author{Abderrahim Mesbah}
\address{Abderrahim Mesbah \newline
Beijing Institute of Mathematical Sciences and Applications, Beijing, China\\
}
\email{abderrahimmesbah@bimsa.cn}

\begin{document}
\maketitle

\begin{abstract}
Let $h^{+}$ and $h^{-}$ be two complete, conformal metrics on the disc $\mathbb{D}$. Assume moreover that the derivatives of the conformal factors of metrics $h^{+}$ and $h^{-}$ are bounded at any order with respect to the hyperbolic metric, and that the metrics have curvatures in the interval $\left(-\frac{1}{\epsilon}, -1 - \epsilon\right)$, for some $\epsilon > 0$. Let $f$ be a quasi-symmetric map. We show the existence of a globally hyperbolic convex subset $\Omega$ (see Definition \ref{globally-hyperbolic-convex-subset}) of the three-dimensional anti-de Sitter space, such that $\Omega$ has $h^{+}$ (respectively $h^{-}$) as the induced metric on its future boundary (respectively on its past boundary) and has a gluing map $\Phi_{\Omega}$ (see Definition \ref{gluing-maps}) equal to $f$.
\end{abstract}

\section{Introduction}
\begin{figure}
    \centering
    \includegraphics[width=18cm]{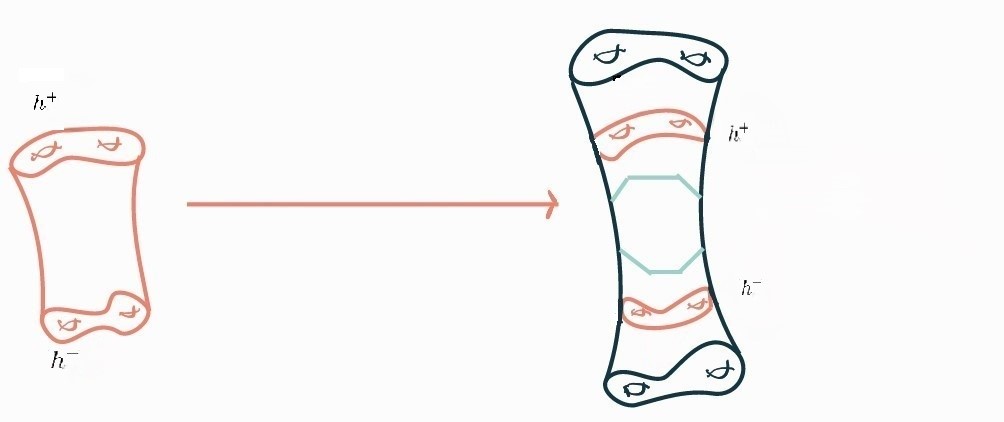}
    \caption{In other words, Tamburelli has shown that for any closed surface $S$ with genus greater than or equal to $2$, and for any Riemannian metrics $h^{+}$ and $h^{-}$ on $S$ with curvatures strictly smaller than $-1$, there exists a maximal globally hyperbolic manifold $M$ homeomorphic to $S \times (0,1)$ and a submanifold $N$ of $M$ that has the same homotopy type as $S$ and is diffeomorphic to $S \times [0,1]$, such that the boundary of $N$ consists of two disjoint space-like surfaces on which $M$ induces a metric on $S \times  \left\{ 0\right\}$ homotopic to $h^{+}$ and a metric on $S \times  \left\{ 1\right\}$ homotopic to $h^{-}$.
}
    \label{existence-pic}
\end{figure}

The three-dimensional anti-de Sitter space $\AdS^{2,1}$ is the Lorentzian analogue of the hyperbolic space $\mathbb{H}^3$. That is, any three-dimensional Lorentzian manifold that has constant sectional curvature equal to $-1$ is locally modeled on $\AdS^{2,1}$. A manifold which is locally modeled on $\AdS^{2,1}$, and which is oriented and time-oriented, is called an $\AdS$ space-time. Particular $\AdS$ space-times, called globally hyperbolic anti-de Sitter manifolds, have attracted attention since the work of \cite{zbMATH05200424} that showed their relation to Teichm\"uller theory and their similarities with quasi-Fuchsian manifolds.\\  
A maximal globally hyperbolic $\AdS$ 3-manifold $M$ is diffeomorphic to $S \times  (0,1)$, where $S$ is a closed connected surface of genus greater or equal to $2$. A maximal globally hyperbolic manifold $\AdS$ 3-manifold is always the quotient of the domain of dependence of a quasi-circle $C_{M}$ (see Definitions \ref{quasi-circle} and \ref{domain-of-dependance}) by a representation of the form $(\rho_1, \rho_2): \pi_1(S) \to \PSL(2,\mathbb{R}) \times \PSL(2,\mathbb{R})$ where $\rho_1$ and $\rho_2$ are Fuchsian representations (see \cite{zbMATH05200424} or \cite[Section 5]{bonsante2020anti}). This gives a one-to-one correspondence between maximal globally hyperbolic $\AdS$ 3-manifolds diffeomorphic to $S \times (0,1)$ and $\mathcal{T}(S) \times \mathcal{T}(S)$, where $\mathcal{T}(S)$ is the Teichm\"uller space of $S$. We have a similar correspondence for quasi-Fuchsian manifolds, which is known as Bers' Theorem (see \cite{bers1960simultaneous}). Also, a globally hyperbolic manifold $M$ contains a smallest non-empty geodesically convex subset $C(M)$ which is called the convex core of $M$. When $C(M)$ is not a totally geodesic two-dimensional manifold, it has a boundary $\partial C(M)$ which consists of the disjoint union of two hyperbolic surfaces pleated along a measured lamination (see \cite{zbMATH05200424} or \cite[Section 6]{bonsante2020anti}).\\
It was proved in \cite{diallo2013prescribing} that given any two hyperbolic metrics on $S$, there exists a globally hyperbolic manifold $M$ that induces these metrics on the boundary $\partial C(M)$. An analogous result is known for quasi-Fuchsian manifolds (it follows from the works of \cite{epstein1987convex} and \cite{labourie1992metriques}). However, in both the globally hyperbolic and quasi-Fuchsian settings, It is not known whether the manifold that induces the given metrics on the boundary of the convex core is unique. Nevertheless, Prosanov proved the uniqueness of the manifold for an open and dense subset of the deformation space of quasi-Fuchsian manifolds (see \cite{Prosanov}).\\
For quasi-Fuchsian manifolds, Thurston conjectured that there is a one-to-one correspondence between the space of quasi-Fuchsian manifolds and the bending laminations of their convex core boundaries. That conjecture has been proved. Indeed, Bonahon and Otal (see \cite{Bonahon-otal}) showed that any two measured laminations on $S$ that fill the surface and do not have a closed leaf with weight greater than or equal to $\pi$ can be realized as the bending laminations of the convex core boundary of some quasi-Fuchsian manifold, which is homeomorphic to $S \times (0,1)$. Later, Dular and Schlenker (see \cite{Dular-schlenker}) showed that if a pair of measured laminations on $S$ arises as the bending laminations of the boundary of the convex core of a quasi-Fuchsian manifold $Q$, then $Q$ is uniquely determined. In other words, the bending laminations uniquely determine the quasi-Fuchsian manifold.\\
For globally hyperbolic $\AdS^{2,1}$ manifolds, a similar result to the one given by Bonahon and Otal was obtained by Bonsante and Schlenker (see \cite{Bonsante-schlenker-realization}). However, it is still unknown whether the bending laminations on the boundary of the convex core uniquely determine the globally hyperbolic manifold, exept if the bending lamination is small enough in a sense given in \cite{Bonsante-schlenker-realization}.
 \\
Later, Tamburelli in \cite{tamburelli2018prescribing} proved that given any two metrics on $S$ that have curvatures stricly smaller than $-1$, we can find a globally hyperbolic $\AdS^{2,1}$ manifold $N$ which is diffeomorphic to $S \times [0,1]$ with smooth space-like boundary (it is embedded in a maximal globally hyperbolic manifold, see Figure \ref{existence-pic}) that induces the given two metrics on $S$ on the boundary of $N$. A similar result was proved by \cite{labourie1992metriques} in the case of quasi-Fuchsian manifolds. Unlike the $\AdS^{2,1}$ case, in which we don't know the uniqueness of the manifold $N$ that induces the given metrics on its boundary, Schlenker in \cite{schlenker2006hyperbolic} showed the uniqueness in the hyperbolic case.\\

The authors in \cite{bonsante2021induced} and \cite{chen2022geometric} have studied lifts of quasi-Fuchsian manifolds and globally hyperbolic manifolds, and related them to the universal Teichm\"uller space.\\
The following theorem is a universal analogue of the main theorem of \cite{diallo2013prescribing} about prescribing the metrics on the boundary of the convex core of a globally hyperbolic manifold. For definition of the gluing map we refer to Definition \ref{gluing-maps}.
\begin{theorem}\cite[Theorem D]{bonsante2021induced}
 Any normalized quasi-symmetric homeomorphism of the circle is realized as the gluing map at infinity for the convex hull of a normalized quasi-circle in $Ein^{1,1} := \partial_{\infty}\AdS^{2,1}$.   
\end{theorem}
The authors in \cite{bonsante2021induced} have also shown a partial analogue of the main theorem in \cite{tamburelli2018prescribing} (when we assume that the curvatures of the induced metrics on the boundary are constant), which prescribes the metrics on the smooth space-like boundaries of a convex compact $\AdS^{2,1}$ globally hyperbolic manifold. In the next theorem we mean by $K$-surface a convex space-like surface (we mean by convex that the product of its principal curvatures is positive) that has constant sectional curvature equal to $K$. 
\begin{theorem}\cite[Theorem E]{bonsante2021induced}
 Any normalized quasisymmetric homeomorphism of the circle is realized as
the gluing map between the future and past $K$ surfaces spanning some normalized quasicircle in $Ein^{1,1} := \partial_{\infty}\AdS^{2,1}$.   
\end{theorem}
The main theorem of this paper extends the main theorem of \cite{tamburelli2018prescribing} to the setting given in \cite{bonsante2021induced}. Our theorem gives a positive answer to the first question of Section 1.7 of \cite{schlenker2020weyl} under the hypothesis that the derivatives of the conformal factors of the metrics $h^{\pm}$ are bounded with respect to the hyperbolic metric.\\
Before stating the theorem, we will define what we mean by saying that the conformal factor of a metric has bounded derivatives with respect to the hyperbolic metric. In this paper, we will simply say that the metric has bounded derivatives, and it should be understood that the derivatives of the conformal factor of the metric are bounded with respect to the hyperbolic metric.
\begin{definition}\label{boundedness}
We denote by $h_{-1}$ the hyperbolic metric on the unit disc $\mathbb{D}$, written in conformal form as $h_{-1} = \frac{4\,|dz|^2}{(1 - |z|^2)^2}$. Let $h = e^{2\rho} h_{-1}$ be a complete conformal metric on $\mathbb{D}$, where $\rho : \mathbb{D} \to \mathbb{R}$ is a smooth function.\\
We say that $h$ has bounded derivatives of order $p$ if there exists a constant $M_p > 0$ such that all derivatives of $\rho$ of order $p$, when measured with respect to the hyperbolic metric $h_{-1}$, are uniformly bounded on the disc by $M_p$, independently of the point $z \in \mathbb{D}$.\\
In terms of partial derivatives, this means that for any $p \in \mathbb{N}$ and for every multi-index $\alpha = (\alpha_1,\alpha_2) \in \mathbb{N}^2$ with $|\alpha| = p$, we have :
$$
\left| D^\alpha \rho(z) \right| \leq M_p (1 - |z|^2)^p \quad \text{for all } z \in \mathbb{D},
$$
where $D^\alpha = \partial_x^{\alpha_1} \partial_y^{\alpha_2}$ and $z = x + iy$. The factor $(1 - |z|^2)^p$ reflects that the derivatives are measured using the norm induced by the hyperbolic metric.\\
We say that $h$ has bounded derivatives if this property holds for every order $p \in \mathbb{N}$.

\end{definition}
Note that any conformal metric on $\mathbb{D}$ that has constant negative sectional curvature satisfies this condition. Also, note that any metric which is invariant under a Fuchsian representation $\rho: \pi_{1}(S) \to \PSL(2,\mathbb{R})$, where $S$ is a closed hyperbolic surface, satisfies this boundedness condition. 
Now we state our main Theorem, for a definition of globally hyperbolic convex subset see Definition \ref{globally-hyperbolic-convex-subset} (see also Figure \ref{Ads-pic}). 
\begin{theorem}\label{mainn}
Let $h^{+}$ and $h^{-}$ be two complete, conformal metrics on the disc $\mathbb{D}$ that have curvatures in an interval of the form $(-\frac{1}{\epsilon},-1-\epsilon)$, for some $\epsilon > 0$. Assume moreover that there is a sequence $(M_p)_{p \in \mathbb{N}}$ of positive real numbers such that any derivative of $h^{+}$ or $h^{-}$ of order $p$ is bounded by $M_{p}$. Let $f$ be a normalized quasi-symmetric map. Then there exists a globally hyperbolic convex subset $\Omega \subset \AdS^{2,1}$ such that the induced metric on $\partial^{+}\Omega$ is isometric to $(\mathbb{D},h^{+})$, the induced metric on $\partial^{-}\Omega$ is isometric to $(\mathbb{D},h^{-})$, and the gluing map of $\Omega$ is equal to $f$.     
\end{theorem}
A natural question is the following
\begin{question}
Is the convex subset $\Omega$ in Theorem \ref{mainn} unique (up to isometry)?
\end{question}
\begin{figure}\label{Ads-pic}
    \centering
    \includegraphics[width=12cm]{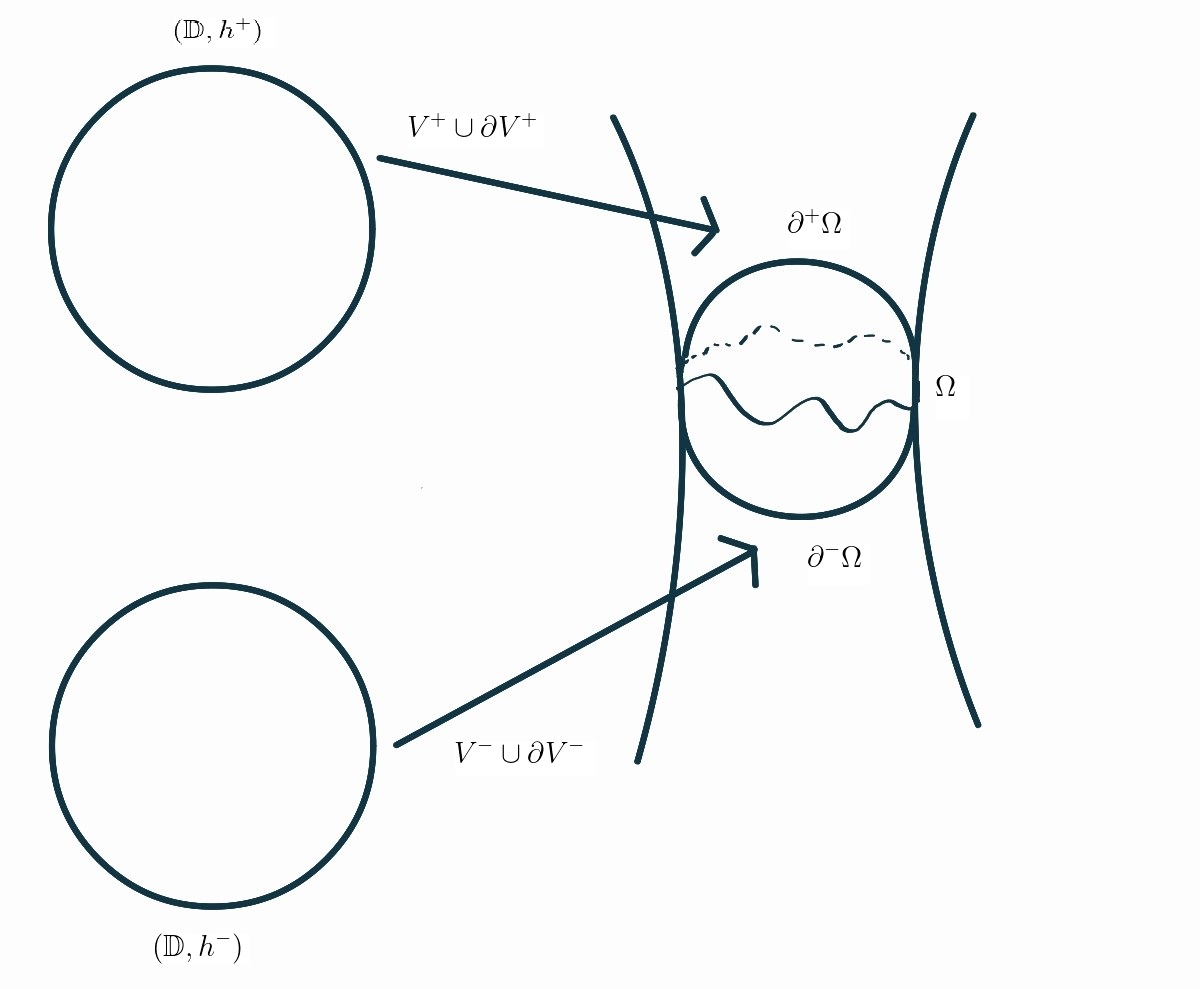}
    \caption{The connected components of $\partial \Omega \cap \AdS^{2,1}$, denoted by $\partial^{+} \Omega$ and $\partial^{-} \Omega$, are isometric to $(\mathbb{D}, h^{+})$ and $(\mathbb{D}, h^{-})$, respectively. The isometries $V^{+}$ and $V^{-}$ can be extended to the ideal boundary and then define the map $(\partial V^{-})^{-1} \circ \partial V^{+}$ (see Section \ref{Gluing-maps-section}), which is referred to as the gluing map of $\Omega$. In Theorem \ref{mainn} this gluing map is equal to $f$.}
    \label{Ads-pic}
\end{figure}

\subsection*{Outline of the paper}
In the first section, we recall the necessary background about quasi-conformal maps, quasi-isometries, quasi-symmetric maps, cross-ratio, and the universal Teichm\"uller space. In the second section, we introduce anti-de Sitter geometry. We present the hyperboloid model, the projective model, and the Lie group model. We discuss concepts such as causality, and meridians. In the third section, we define globally hyperbolic anti-de Sitter 3-manifolds and introduce the notion of globally hyperbolic convex subsets of $\AdS^{2,1}$. In the fourth section, we study the geometry of the boundary of these globally hyperbolic convex subsets, showing that the asymptotic behavior near the boundary at infinity defines a quasi-symmetric map, referred to as the gluing map. In the fifth section, we approximate any metric on the disc with bounded negative curvature, where its derivatives are bounded with respect to the hyperbolic metric, by metrics invariant under Fuchsian representations that have uniformly bounded curvatures and uniformly bounded derivatives. In the final section, we prove the main theorem by using an approximation involving lifts of globally hyperbolic $\AdS$ 3-manifolds.

\subsection*{Acknowledgments}
I would like to express my sincere gratitude to my supervisor, Jean-Marc Schlenker, for his invaluable assistance and unwavering support throughout this work.\\
Additionally, I am deeply thankful to Francesco Bonsante for his insightful guidance during my stay in Pavia. I am profoundly grateful to Nathaniel Sagman for providing the proof of Lemma \ref{nathanial} and for his valuable remarks and comments.\\
This work has been supported by the Luxembourg National Research Fund PRIDE/17/1224660/GPS.

\section{Preliminaries}
\subsection{Quasi-conformal maps}\label{1.6}
Let $X$ and $Y$ be Riemann surfaces (not necessarily compact). Let\\ $f: X \to Y$ be an orientation preserving diffeomorphism. We define the Beltrami differential $\mu = \mu(f)$ by the equation $\frac{\partial f}{\partial \Bar{z}} = \mu \frac{\partial f}{\partial z}$. We say that $f$ is $K$ quasi-conformal if the dilatation number $K(f) = \frac{1 + \left| \mu \right|_{\infty}}{1 - \left| \mu \right|_{\infty}}$ is less than or equal to $K$. Note that we don't need $f$ to be a $C^{1}$ diffeomorphism to define the notion of quasi-conformal maps. In fact, all we need is for $f$ to be a homeomorphism between $X$ and $Y$ that has derivatives in the sense of distribution that are $L^{2}$. For more details see \cite{lehto1973quasiconformal}.
The following proposition is well known, see for example \cite{teleman2007fletcher}.
\begin{proposition}\label{qfextends}
Any quasi-conformal homeomorphism $f: \mathbb{H}^2 \to \mathbb{H}^2$ extends to a homeomorphism $\partial f : \partial_{\infty}\mathbb{H}^2 \to \partial_{\infty}\mathbb{H}^2$.    
\end{proposition}
\subsection{Quasi-isometries}\label{1.7}
Let $(X,d_{X})$ and $(Y,d_{Y})$ be two metric spaces. Let $A \geq 1$ and $B \geq 0$. We say that a map $f: (X,d_{X}) \to (Y,d_{Y})$ is a $(A,B)$ quasi-isometric embedding if for any $x_1$ and $x_2$ in $X$ the following inequalities hold:
$$\frac{1}{A}d_{X}(x_1,x_2) - B \leq d_{Y}(f(x_1),f(x_2)) \leq Ad_{X}(x_1,x_2) + B.  $$
Let $C \geq 0$. We say that $f$ is an $(A,B,C)$ quasi-isometry if it is an $(A,B)$ quasi-isometric embedding and it is $C$-dense, that is:
$$\forall z \in Y, \exists \ x \in X, d_{Y}(f(x),z) \leq C  $$
If $(X,d_{X})$ and $(Y,d_{Y})$ are $\delta$-hyperbolic spaces, then any quasi-isometric embedding $f$ extends uniquely to a homeomorphism $\partial f: \partial_{\infty}X \to \partial_{\infty}Y$ of the visual boundary. The homeomorphism $\partial f$ is called a quasi-symmetric map and has many interesting properties. We will discuss the notion of quasi-symmetric maps next.
It is worth to mention the following well known proposition, for a proof see for example \cite{teleman2007fletcher}.
\begin{proposition}\label{qf-to-qi}
 Any quasi-conformal map $f: \mathbb{H}^2 \to \mathbb{H}^2$ is a quasi-isometric embedding.    
\end{proposition}
Note that Proposition \ref{qf-to-qi} implies Proposition \ref{qfextends}.

\subsection{Quasi-symmetric maps}\label{1.8}
For more details see \cite{hubbard2016teichmuller}. We denote $\mathbb{RP}^1 := \mathbb{R} \cup \left\{\infty \right\}$. Let $\phi: \mathbb{RP}^{1}  \to \mathbb{RP}^{1}  $ be a strictly increasing homeomorphism that satisfies $\phi(\infty) = \infty$. We say that $\phi$ is quasi-symmetric if there exists $k > 0$ such that,
$$\forall x \in \mathbb{R}, \ \forall t \in \mathbb{R}^{*}_{+}, \ \frac{1}{k} \leq\frac{\phi(x+t)-\phi(x)}{\phi(x) - \phi(x-t)} \leq k .  $$
In this case we say that $\phi$ is $k$ is quasi-symmetric, and we call $k$ by the quasi-symmetric constant of $\phi$.\\
If $\phi$ does not fix $\infty$, then we say that $\phi$ is $k$ quasi-symmetric if there exists an element $g \in PSL(2,\mathbb{R})$ (therefore many elements) such that $g \circ \phi (\infty) = \infty$ and $g \circ \phi$ is $k$ quasi-symmetric.\\
The quasi-symmetric maps can be seen as the extension of quasi-conformal maps as shows the following proposition.
\begin{proposition}\label{hubbard}
\noindent
\begin{itemize}
    \item For any $k \geq 1 $ there exists $k' > 0$ such that any $K$ quasi-conformal map $f: \mathbb{H}^2 \to \mathbb{H}^2$ has a continuous extension to the boundary $\partial f : \mathbb{RP}^1  \to \mathbb{RP}^1 $ which is $k'$ quasi-symmetric.
    \item For any $M > 0$ there exists $M' > 1$ such that any $k'$ quasi-symmetric map $\phi: \mathbb{RP}^1  \to \mathbb{RP}^1$ has a continuous extension $f: \mathbb{H}^2 \to \mathbb{H}^2$ which is $M'$ quasi-conformal.
\end{itemize}
\end{proposition}
    For the first point see for example \cite[Corollary 4.9.4]{hubbard2016teichmuller}.\\
    For the second point see for example The Douady-Earle extension theorem, \cite[Theorem 5.1.2]{hubbard2016teichmuller}   

Also, quasi-symmetric maps can be seen as the extension of quasi-isometric embeddings as shows the following proposition. 
\begin{proposition}
a map $f : \partial_{\infty}\mathbb{H}^2 \to \partial_{\infty}\mathbb{H}^2$ is quasi-symmetric if and only if there exists a quasi-isometric embedding $F: \mathbb{H}^2 \to \mathbb{H}^2$ such that $f = \partial F$ 
\end{proposition}
Since any quasi-conformal map is a quasi-isometry then we have one implication.\\
For the second implication we refer for example to \cite{hubbard2016teichmuller} or \cite{teleman2007fletcher}.
    
Fix $k >0$, an interesting property of $k$ quasi-symmetric maps is the compactness as show the following proposition, see for example \cite[Corollary 4.9.7]{hubbard2016teichmuller}.
\begin{proposition}\label{compact}
  Let $(\phi_{n})_{n \in \mathbb{N}}$ be a sequence of $k$ quasi-symmetric maps such that, for any $n$, $\phi_{n}(i) = i$, for any $i \in \left\{ 0,1,\infty \right\}$ (We say that such a map is normalized). Then up to extract a subsequence, $(\phi_{n})_{n \in \mathbb{N}}$ converge  to a $k$ quasi-symmetric map $\phi_{\infty}$ in the topology $C^{0}(\mathbb{RP}^1)$.    
\end{proposition}
\subsection{Cross-ratio}\label{1.9}
The cross-ratio of four points $(a,b,c,d) \in (\partial_{\infty}\mathbb{H}^2)^{4}$ is defined by the formula (note that the order of $(a,b,c,d)$ is important):
$$cr(a,b,c,d) = \frac{(c-a)(d-b)}{(b-a)(d-c)}.$$
It is well known that the cross-ratio is invariant under the action of $PSL(2,\mathbb{R})$, that is if $g \in PSL(2,\mathbb{R})$, then $cr(a,b,c,d) = cr(g(a),g(b),g(c),g(d))$.\\
A quadruple of points $(a,b,c,d)$ is called symmetric if $cr(a,b,c,d) = -1$, or equivalently if there is $g \in PSL(2,\mathbb{R})$ such that $(g(a),g(b),g(c),g(d)) = (0,1,-1,\infty)$.\\
For an orientation preserving homeomorphism $f: \partial_{\infty} \mathbb{H}^2 \to \partial_{\infty} \mathbb{H}^2$, we denote:  
$$crf(a,b,c,d) := cr(f(a), f(b), f(c), f(d)).$$
The following characterisation of quasi-symmetric maps using cross-ratio is well known, see \cite{teleman2007fletcher} for example.
\begin{theorem}
For any $k \geq 1$, there exists $M \geq 1$ such that if $f: \partial_{\infty}\mathbb{H}^2 \to \partial_{\infty}\mathbb{H}^2$ is $k$ quasi-symmetric, then for any symmetric quadruple $(a,b,c,d)$ the following inequalities hold:
$$-M \leq crf(a,b,c,d) \leq -\frac{1}{M}.$$
Moreover, if $k$ goes to infinity, then so does $M$.\\
Conversely for any $M \geq 1$ there exists $k \geq 1$ such that if $f: \partial_{\infty}\mathbb{H}^2 \to \partial_{\infty}\mathbb{H}^2$ is an orientation preserving homeomorphism that satisfies :
$$-M \leq crf(a,b,c,d) \leq -\frac{1}{M},$$
then $f$ is $k$ quasi-symmetric.
Again, if $M$ goes to infinity, then so does $k$.
\end{theorem}
\subsection{The universal Teichm\"{u}ller space}\label{1.10}
The quasi-symmetric maps of $\partial_{\infty}\mathbb{H}^2$ form a group that we denote by $\mathcal{QS}(\partial_{\infty}\mathbb{H}^2)$. Recall that $PSL(2,\mathbb{R})$ acts on $\mathcal{QS}(\partial_{\infty}\mathbb{H}^2)$ by post-composition. We define the universal Teichm\"uller space to be:
$$\mathcal{T} = \mathcal{QS}(\partial_{\infty}\mathbb{H}^2) / PSL(2,\mathbb{R}).$$
Or equivalently, we can identify $\mathcal{T}$ with the set of normalized quasi-symmetric maps, that is, the set of quasi-symmetric maps $f$ such that $f(p) = p$ for any $ p \in \left\{0,1,\infty \right\}$.\\
For any genus $g \geq 2$, the classical Teichm\"uller space $\mathcal{T}(S)$ is embedded in the universal Teichm\"uller space. Indeed, let us fix one Riemann structure $X$ on $S$ (which serves as a base point for $\mathcal{T}(S)$). $\tilde{X}$, the universal cover of $X$, is conformal to $\mathbb{H}^2$, so we can simply identify it with $\mathbb{H}^2$. This implies that we can identify $X$ with $\mathbb{H}^2 / \Gamma$.\\
Suppose $Y$ is another Riemann structure on $S$, and let $f: X \to Y$ be a diffeomorphism. We can identify $X$ with $\mathbb{H}^2 / \Gamma$ and $Y$ with $\mathbb{H}^2 / \Gamma'$, where $\Gamma$ and $\Gamma'$ are Fuchsian groups. Let $\tilde{f}: \mathbb{H}^2 \to \mathbb{H}^2$ be a lift of $f$. Since $S$ is compact, $f$ is quasi-conformal. Hence, so is $\tilde{f}$, by Proposition \ref{hubbard}. $\tilde{f}$ extends to a unique quasi-symmetric map $\partial \tilde{f}$. Note that adjusting $X$ or $Y$ by isotopy changes $\partial \tilde{f}$ up to M\"obius translation. In conclusion, once we fix a base point $X$ in $\mathcal{T}(S)$, we obtain a unique (up to composition by M\"obius maps) quasi-symmetric map for each $Y \in \mathcal{T}(S)$. This correspondence is an embedding because for each two Fuchsian representations $\rho_{1}, \rho_{2}: \pi_{1}(S) \to PSL(2,\mathbb{R})$, there exists a unique quasi-symmetric map which is $\rho_{1}$-$\rho_{2}$ equivariant (that is, for all $\gamma \in \pi_1(S)$ and all $x \in \partial_{\infty} \mathbb{H}^2$, we have $\rho_1(\gamma) \cdot f(x) = f(\rho_2(\gamma) \cdot x)$). As a consequence, each point $Y$ in $\mathcal{T}(S)$ corresponds to a unique normalized quasi-symmetric map.
\section{Anti-de Sitter space}
The anti-de Sitter space $\AdS^{n,1}$ is the Lorentzian analogue of the hyperbolic space $\mathbb{H}^{n+1}$, in the sense that any Lorentzian manifold of constant sectional curvature equal to $-1$ is locally modeled on $\AdS^{n,1}$.
In this section we will introduce the models for the three dimensional anti-de Sitter space $\AdS^{2,1}$ that we will use later thought this paper. For more details we refer the reader to \cite{bonsante2020anti} and\cite{zbMATH05200424}.
\subsection{Hyperbloid model}
Let $q_{2,2}$ be the quadratic form defined on $\mathbb{R}^{4}$ by the formula :
$$q_{2,2}(x_1,x_2,x_3,x_4) = x^2_1 + x^2_2 - x^2_3 - x^2_4.$$
We denote by $\mathbb{R}^{2,2}$ the space $(\mathbb{R}^{4},q_{2,2})$.\\
We define the space $\mathbb{H}^{2,1}$ to be :
$$\mathbb{H}^{2,1} :=  \left\{x \in \mathbb{R}^{4}, q_{2,2}(x) = -1 \right\}.$$
Then $q_{2,2}$ induce a scalar product on each tangent space of $\mathbb{H}^{2,1}$ that has signature $(2,1)$, this makes $\mathbb{H}^{2,1}$ a Lorentian manifold. We refer to \cite{bonsante2020anti} and\cite{zbMATH05200424} to see why $\mathbb{H}^{2,1}$ has a constant sectional curvature equal to $-1$.
\subsection{The projective model of $\AdS^{2,1}$}
We introduce $\AdS^{2,1}$, the projective model of the anti-de Sitter space, to be :
$$\AdS^{2,1} := \mathbb{H}^{2,1} /  \left\{\pm \right\},$$
or equivalently :
$$\AdS^{2,1} := \left\{\left [ x \right ] \in \mathbb{RP}^3, q_{2,2}(x) < 0  \right\}.$$
The projective model allows us to visualise better $\partial_{\infty} \AdS^{2,1}$, the ideal boundary of the anti-de Sitter space.
$$\partial_{\infty} \AdS^{2,1} := \left\{\left [ x \right ] \in \mathbb{RP}^3, q_{2,2}(x) = 0  \right\}$$
The anti-de Sitter space $\AdS^{2,1}$ induces a Lorentzian conformal structure on its ideal boundary $\partial_{\infty}\AdS^{2,1}$ (see \cite[Section 2.2]{bonsante2020anti}).
\subsection{The Lie group model}
Let $\mathcal{M}_{2,2}(\mathbb{R})$ be the space of $2 \times 2$ matrices with real coefficient. The space $(\mathcal{M}_{2,2}(\mathbb{R}), -det)$ is isometric to $\mathbb{R}^{2,2}$ via the map
\begin{align*}
 \mathbb{R}^{4} & \to \mathcal{M}_{2,2}(\mathbb{R})\\
 (x_1,x_2,x_3,x_4) & \mapsto \begin{pmatrix}
x_1 - x_3 & x_4 - x_2  \\
x_2 + x_4 & x_1 + x_3  \\
\end{pmatrix}     
\end{align*}

and under this isomorphism $\mathbb{H}^{2,1}$ is identified with the Lie group $SL(2,\mathbb{R})$ (see for example \cite[Section 2.1]{bonsante2021induced}). Under this identification it yields that the projective model $\AdS^{2,1}$ is identified with $\PSL(2,\mathbb{R})$ and $\partial_{\infty}\AdS^{2,1}$ is identified with $\left\{\left [ M \right ] \in \PSL_{2}(\mathbb{R}), det(M) = 0  \right\}$.\\
There is an explicite identification between $\partial_{\infty}\AdS^{2,1}$ and $\mathbb{RP}^{1} \times \mathbb{RP}^1$ via the following map :
\begin{align*}
 \partial_{\infty} \AdS^{2,1} & \to \mathbb{RP}^1 \times \mathbb{RP}^1\\
 \left [ M \right ] &  \to (Im(M),Ker(M)).
\end{align*}
Note that group $\PSL(2,\mathbb{R}) \times \PSL(2,\mathbb{R})$ acts on $\PSL(2,\mathbb{R})$ by left and right composition, that is :
$$(A,B).X = AXB^{-1}.$$
We refer the reader to \cite[Section 3.1]{bonsante2020anti} to see why we can identify the isometry group of $\AdS^{2,1}$ that preserve orientation and time orientation with $\PSL(2,\mathbb{R}) \times \PSL(2,\mathbb{R})$.\\
In general, in any Lorentzian manifold $(M,q)$, we say that a vector $v \in T_{p}M$ is :
\begin{itemize}
    \item Space-like if $q(v,v) > 0$.
    \item Light-like if $q(v,v) = 0$.
    \item Time-like if $q(v,v) < 0$.
\end{itemize}
The geodesics in $\AdS^{2,1}$ are obtained by the intersection of planes of $\mathbb{R}^{2,2}$ that go thought the origin with $\mathbb{H}^{2,1}$. We say that a geodesic $\alpha$ is :
\begin{itemize}
    \item Space-like if $q_{2,2}(\Dot{\alpha}) > 0$.
    \item Light-like if $q_{2,2}(\Dot{\alpha}) = 0$.
    \item Time-like if $q_{2,2}(\Dot{\alpha}) < 0$.
\end{itemize} 
We refer to \cite[Section 2.3]{bonsante2020anti} to see that $\AdS^{2,1}$ is time oriented.
\subsection{Totally geodesic planes}
A plane $P$ in $\AdS^{2,1}$ is defined as the intersection of $W$, a three dimensional vector space in $\mathbb{R}^{2,2}$, with $\mathbb{H}^{2,1}$. That is a plane $P$ in $\AdS^{2,1}$ is defined as : 
$$P := W \cap \mathbb{H}^{2,1}.$$
We say that a plane $P$ is :
\begin{itemize}
    \item Space-like if any two points in $P$ are connected by a space-like geodesic of $\AdS^{2,1}$ which is in $P$.
    \item Time-like if it contains a time-like geodesic.
    \item Light-like otherwise.
\end{itemize}

\subsection{Acausal curves}
Let $C$ be a continuous curve in $\partial^{\infty} \AdS^{2,1}$.
We say that $C$ is achronal (resp acausal), if for any point $p \in C$, there is a neighborhood $U$ of $p \in \partial_{\infty} \AdS^{2,1}$, such that $U \cap C$ is contained in the complement of the regions of $U$ which are connected to $p$ by timelike curves (resp timelike and lightlike curves).\\
We have seen that $\partial_{\infty}\AdS^{2,1}$ is identified with $\mathbb{RP}^1 \times \mathbb{RP}^1$. Then the graph of any homeomorphism $f: \mathbb{RP}^1 \to  \mathbb{RP}^1$ defines a curve on $\partial_{\infty}\AdS^{2,1}$.
\begin{definition}\label{quasi-circle}
An acausal curve $C \subset \partial_{\infty}\AdS^{2,1}$ is a quasi-circle if it is the graph of a quasi-symmetric map.    
\end{definition}\label{domain-of-dependance}
We will define the domain of dependence of a quasi-circle $C$. We say that a curve in $\AdS^{2,1}$ is causal if its tangent vector at any point is time-like or light-like.
\begin{definition}
Let $C \subset \partial_{\infty}\AdS^{2,1}$ be a quasi-circle. We define $D(C)$, the domain of dependence of $C$, to be:
$$
D(C) := \left\{\, p \in \AdS^{2,1} \;\middle|\; p \text{ is connected to } C \text{ by no causal path} \,\right\}.
$$

Equivalently, $D(C)$ (see \cite{BB09}) is the unique maximal (in the sense of inclusion) open convex subset whose boundary at infinity is equal to $C$. 
\end{definition}
For any quasi-circle $C$, the domain of dependence $D(C)$ is contained in a unique affine chart and it is the maximal convex subset of $\AdS^{2,1} \cup \partial_{\infty}\AdS^{2,1}$ that contains $C$.\\
We also give the following definition.
\begin{definition}
  Let $C \subset \partial_{\infty}\AdS^{2,1}$ be a quasi-circle. We define $\CH(C)$ to be the smallest closed convex subset that contains $C$.  
\end{definition}
If $C$ is not the graph of a projective map then $\partial \CH(C) \cap \AdS^{2,1} = \partial^{+}\CH(C) \bigsqcup \partial^{-}\CH(C)$, where each of $\partial^{\pm} CH(C)$ is a topological disc which is a pleated surface.

\subsection{The width of the convex hull of a meridian}
The width of a quasi-circle is the defined to be the maximum time like distance between the boundaries of its convex hull.
\begin{definition}
  Let $C \subset \partial_{\infty} \AdS^{2,1}$ be an achronal meridian. The width $w(C)$ of $C$ is the supremum of the time distance between a point of $\partial^{-}CH(C)$ and $\partial^{+}CH(C)$.
\end{definition}
In the next proposition, we will use the notion of a rhombus in $\partial_{\infty}\AdS^{2,1}$, so we define it below before stating the proposition. See also \cite[Example 6.7]{bonsante2021induced}.
\begin{definition}\label{rhombus}
Recall that $\partial_{\infty} \AdS^{2,1}$ is identified with $\mathbb{RP}^1 \times \mathbb{RP}^1$. We call the lines of the form $\{p\} \times \mathbb{RP}^1$ (resp $\mathbb{RP}^1 \times \{p\}$) the right (resp left) ruling lines, where $p \in \mathbb{RP}^1$.\\
Let $a, b, a', b' \in \mathbb{RP}^1$ such that $a < a'$ and $b < b'$ in the cyclic order of $\mathbb{RP}^1$. We denote by $C_{\diamond}$ the piecewise linear curve in $\partial_{\infty} \AdS^{2,1}$ consisting of the horizontal segments (that is, segments lying on a right ruling line) connecting $(a,b)$ to $(a,b')$ and $(a',b)$ to $(a',b')$, and the vertical segments (that is, segments lying on a left ruling line) connecting $(a,b)$ to $(a',b)$ and $(a,b')$ to $(a',b')$. We call such a curve a \emph{rhombus}.\\ 
Note that there is a unique rhombus up to isometry (that is, up to the action of isometries of $\AdS^{2,1}$). The convex hull $\CH(C_{\diamond})$ of a rhombus in $\partial_{\infty} \AdS^{2,1}$, taken inside $\AdS^{2,1}$, is a tetrahedron (see \cite[Figure 3]{bonsante2021induced}). Since the spacelike segment (inside $\AdS^{2,1}$) connecting $(a,b)$ to $(a',b')$ is dual to the spacelike segment connecting $(a',b)$ to $(a,b')$, it follows that $w(C_{\diamond}) = \frac{\pi}{2}$. Therefore, the width of a rhombus curve in $\partial_{\infty} \AdS^{2,1}$ is equal to $\frac{\pi}{2}$.
\end{definition}

\begin{proposition}\cite{bonsante2010maximal}\label{rhombus}
 Let $C \subset \partial_{\infty} \AdS^{2,1}$ be an acausal curve. Then $C$ is a quasicircle if and only if $w(C) < \frac{\pi}{2}$. Further, if $C_{n}$ is a sequence of quasicircles whose optimal quasisymmetric constant diverges to infinity, then there exist isometries $\phi_{n} \in Isom({\AdS^{2,1}})$ so that $\phi_{n}(C_n)$ converges to the rhombus $C^{*}$, so that in particular $w(C_n) \to \frac{\pi}{2}$.
   
\end{proposition}
\subsection{The left and right projections}
We begin the subsection by defining a spacelike surface.
\begin{definition}
Let $S$ be a surface and let $(M,g)$ be a Lorentzian manifold, a $C^{1}$ immersion $\sigma : S \to M$ is said to be spacelike if the pull back metric $\sigma^{*}(g)$ is a Riemannian metric. Moreover an immersed surface is convex if its  principal curvatures have the same sign (the sign can be positive or negative depending on the orientation).    
\end{definition}
Let $S$ be a convex space like surface in $\AdS^{2,1}$. Here we introduce the left and right projections $\Pi_{l}$ and $\Pi_{r}$ associated to $S$.\\
Recall (see \cite[Proposition 3.5.2]{bonsante2020anti} for example) that any time-like geodesic $L$ in $\AdS^{3}$ is of the form $L = \left\{A \in PSL(2,\mathbb{R}), Ax = x' \right\}$ for some unique $x,x' \in \mathbb{H}^{2}$. Since at each point $p$ of the surface $S$ there is a unique time-like geoedsic $L_{p}$ that goes thought $p$ and orthogonal to $T_{p}S$, we have a well defined map :
\begin{align*}
  \Pi : S \to & \mathbb{H}^2 \times \mathbb{H}^2 \\ 
  p \mapsto &  (\Pi_{l}(p),\Pi_{r}(p)).
\end{align*}
Where $\Pi_{l}(p),\Pi_{r}(p)$ are defined as $L_{p} =  \left\{A \in PSL(2,\mathbb{C}), A\Pi_{l}(p) = \Pi_{r}(p) \right\}$. It was shown that when the surface $S$ is locally convex, the map $\Pi: S \to \mathbb{H}^{2} \times \mathbb{H}^{2}$ is a local diffeomorphism. Moreover Krasnov and Schlenker in \cite{krasnov2007minimal} computed the pullback of the hyperbolic metric thought $\Pi_{l}$ and $\Pi_{r}$ in terms of the embedding data of the immersion, denote the hyperbolic metric by $h_{-1}$.
\begin{theorem}\cite{krasnov2007minimal}\label{estimation-krasnov-schlenker}
The following equalities hold,
 \begin{align*}
\Pi^*_{l}(h_{-1})(v, w) &= I\left((E + J_{I}B)v, (E + J_{I}B)w\right), \\
\Pi^*_{r}(h_{-1})(v, w) &= I\left((E - J_{I}B)v, (E - J_{I}B)w\right),
\end{align*}  
where $E$ denotes the identity operator, and $J_{I}$ is the complex structure over $TS$ induced by $I$. (Note that while both $J_{I}$ and $B$ depend on the choice of an orientation on $S$, the product is independent of the orientation).
\end{theorem}
 
\begin{lemma}\label{proj extetion}\cite[Lemma 3.18 and Remark 3.19]{bonsante2010maximal}
Assume that $S \subset \AdS^{2,1}$ is a properly embedded spacelike convex surface whose accumulation set in $\partial_{\infty} \AdS^{2,1} = \text{Ein}^{1,1} = \mathbb{RP}^1 \times \mathbb{RP}^1$ is an acausal meridian $C = \Gamma(f)$. Then the maps $\Pi_l$ and $\Pi_r$ are diffeomorphisms onto $\mathbb{H}^2$ that continuously extend the canonical projections $\pi_l, \pi_r : C \to \mathbb{RP}^1$ defined by $\pi_l(x, f(x)) = x$ and $\pi_r(x, f(x)) = f(x)$.
\end{lemma}
\section{Globally hyperbolic convex subsets}
We begin this section by giving the definition of globally hyperbolic manifolds.
\begin{definition}\label{globally-hyperbolic-convex-subset}
Let $M$ be a Lorentzian manifold and let $X$ be an achronal subset of $M$. We say that $X$ is a Cauchy surface if every inextensible causal curve of $M$ meets $X$. We say that a Lorentzian manifold $M$ is globally hyperbolic if it contains a Cauchy surface $X$.
\end{definition}
We refer to \cite{beem1981global}, \cite{geroch1970domain}, \cite{zbMATH02212120}, and \cite{zbMATH02212120} for details for the next theorem.
\begin{theorem}
Let $M$ be a globally hyperbolic manifold, and let $X$ be a Cauchy surface. Then:
\begin{itemize}
    \item Any two Cauchy surfaces of $M$ are diffeomorphic.
    \item There exists a submersion $\phi : M \to \mathbb{R}$ whose fibers are Cauchy surfaces.
    \item $M$ is diffeomorphic to $X \times \mathbb{R}$.
\end{itemize}
\end{theorem}
Next we define the notion of a globally hyperbolic convex subset : 
\begin{definition}
 Let $\Omega$ be a convex subset of $\AdS^{2,1} \cup \partial_{\infty} \AdS^{2,1}$ such that
 \begin{itemize}
     \item $\Omega$ is homemorphic to the ball.
     \item $\partial_{\infty} \Omega$ is a quasi-circle.
     \item $\partial \Omega \cap \AdS^{2,1}$ is the disjoint union of two smooth spacelike disks $\partial^{\pm} \Omega$.
     \item The induced metric on each of $\partial^{\pm} \Omega$ is complete and has curvature in the interval $(-\frac{1}{\epsilon}, -1-\epsilon)$.
 \end{itemize}
 We call a convex $\Omega$ with the preceding properties by a \textbf{globally hyperbolic convex subset}.
\end{definition}
Let $S$ be a closed hyperbolic surface. Then any convex globally hyperbolic manifold $S \times [0,1]$ with a smooth boundary is the quotient of a globally hyperbolic convex subset.

\begin{theorem}\cite{zbMATH05200424}
Let $S \times \left [ 0,1 \right ]$ be an $\AdS^{2,1}$ convex globally hyperbolic manifold. Then there are two Fuchsian representations $\rho_1, \rho_2 : \pi_1(S) \to \PSL(2,\mathbb{R})$ and a globally hyperbolic convex subset $\Omega_{\rho_1, \rho_2}$ such that $S \times \left [ 0,1 \right ]$ is isometric to the quotient $\Omega_{\rho_1, \rho_2} / (\rho_1, \rho_2)$.
\end{theorem}

The author in \cite{tamburelli2018prescribing} has proven the following theorem.

\begin{theorem}\cite{tamburelli2018prescribing}\label{tumb}
Given two metrics $h^+$ and $h^-$ with curvatures $\kappa^{\pm} < -1$ on a closed, oriented surface $S$ of genus $g \geq 2$, there exists an $\AdS^{2,1}$ convex globally hyperbolic manifold $S \times \left [ 0,1 \right ]$ with a smooth, spacelike, strictly convex boundary such that the induced metric on $S \times \{0\}$ is isotopic to $h^{-}$ and the induced metric on $S \times \{1\}$ is isotopic to $h^{+}$.
\end{theorem}

\section{Gluing maps}\label{Gluing-maps-section}
\subsection{Definition of gluing maps}
Let $\Omega$ be a globally hyperbolic convex subset. The induced metric on the boundary of $\Omega$ gives two metrics on the disc $\mathbb{D}$, and the asymptotic behavior of these two metrics near the ideal boundary of $\Omega$ gives a quasi-symmetric map that we will call the gluing map. In this section, we will show the existence and define the gluing maps.\\
Under some conditions on the induced metrics on the boundary of $\Omega$, the principal curvatures of $\partial^{\pm} \Omega$ are always bounded.\\
Before stating the next lemma and proposition, we draw the reader's attention to the fact that in the $\PSL(2,\mathbb{R})$ model of $\AdS^{2,1}$, a lightlike plane $P$ is the orthogonal to an element $\left\lfloor A \right\rfloor \in \PSL(2,\mathbb{R})$ such that $\det(A) = 0$ (see \cite[Section 3.5]{diaf-seppi}). Moreover, the boundary at infinity of a lightlike plane is given by 
$\partial_{\infty}P := (\operatorname{Im}(A) \times \mathbb{RP}^1) \cup (\mathbb{RP}^1 \times \operatorname{Ker}(A))$ (see \cite[Lemma 3.5]{diaf-seppi}).
\begin{lemma}\cite[Lemma 6.3]{bonsante2021induced}\label{compactness-of-quasi-circle}
 Let $C_n$ be a sequence of $k$-quasicircles in $\partial_{\infty}\AdS = \mathbb{RP}^1 \times \mathbb{RP}^1$.\\  
Then there exists a subsequence which converges in the Hausdorff topology to either a  $k$-quasi circle $C$, or to the union $ \left\{p \right\} \times \mathbb{RP}^1 \cup \mathbb{RP}^1 \times  \left\{ q\right\}$ of a line of the left ruling and a line of the right ruling (see Definition \ref{rhombus} for definition of left ruling and right ruling lines).
   
\end{lemma}
\begin{proposition}\cite[Lemma 7.9]{bonsante2021induced}\label{Hausdorff-converge}
Let $S_{n}$ be a sequence of properly embedded convex spacelike disks spanning a sequence of $k$-quasicircles $C_n$. If $C_n$ converges in the Hausdorff topology to the union $\{p\} \times \mathbb{RP}^1 \cup \mathbb{RP}^1 \times \{q\}$ of a line of the left ruling and a line of the right ruling, then $S_n$ converges in the Hausdorff topology to the lightlike plane with boundary at infinity $\{p\} \times \mathbb{RP}^1 \cup \mathbb{RP}^1 \times \{q\}$. If $C_n$ converges in the Hausdorff topology to a $k$-quasicircle $C$, then, up to extracting a subsequence, $S_n$ converges in the Hausdorff topology to a locally convex properly embedded surface spanning the curve $C$.
\end{proposition}
In what follows, we need Theorem 5.5 and Theorem 5.6 from \cite[Section 5]{schlenker1996surfaces}. They are not stated in their full generality, we have stated only what is needed for the rest of our proof. Refer to the same reference for more definitions and details. In particular, for us, $(M, g_0)$ is assumed to be $\AdS^{2,1}$, but in \cite{schlenker1996surfaces}, it is assumed to be a more general $3$-dimensional Lorentzian manifold.
\begin{theorem}\label{JM-thm5.6}\cite[Theorem 5.5]{schlenker1996surfaces}
 Let $f_n : \mathbb{D} \rightarrow M$ be a sequence of uniformly elliptic space-like immersions of the disk such that $(f_n^*( g_0))_{n \in \mathbb{N}}$ converges in $C^\infty$ to a metric $g_\infty$, and such that the sequence of integrals of the mean curvatures $H_n$ is bounded by a constant. Suppose also that $(x_n)_{n \in \mathbb{N}}$ and $(j^{1}f_n(x_n))_{n \in \mathbb{N}}$ converge. Then there exists a subsequence of $(f_n)_{n \in \mathbb{N}}$ that converges in $C^\infty$ to an isometric immersion $f_{\infty}$.   
\end{theorem}
Later, Schlenker studies what happens when the assumptions of Theorem \ref{JM-thm5.6} do not hold.
\begin{theorem}\label{JM-nothm5.6}\cite[Theorem 5.6]{schlenker1996surfaces}
 Let $f_n : \mathbb{D} \to M$ be a sequence of uniformly elliptic space-like immersions of the disk, such that $(f_n^*( g_0))_{n \in \mathbb{N}}$ converges in $C^\infty$ to a metric $g_\infty$. Let $(x_n)_{n \in \mathbb{N}}$ be a sequence converging to $x_{\infty}$ such that $(j^{1}f_n(x_n))_{n \in \mathbb{N}}$ converges, but $(f_n)_{n \in \mathbb{N}}$ does not converge in $C^\infty$ in a neighborhood of $x_\infty$. Then there exists a subsequence of $(f_n)_{n \in \mathbb{N}}$, still denoted by $(f_n)_{n \in \mathbb{N}}$, a maximal geodesic $\gamma$ through $x_\infty$ in $(\mathbb{D},g_{\infty})$, and a geodesic segment $\Gamma$ in $M$, such that $(f_n|\gamma)$ converges to an isometry onto $\Gamma$.   
\end{theorem}

Then we use similar arguments to \cite[Proposition 7.10]{bonsante2021induced} and \cite[Lemma 7.11]{bonsante2021induced} to show the following proposition.\\

Before proceeding with the proof, we briefly recall the notion of the third fundamental form.\\
Let $S$ be a smooth, space-like embedded surface in $\AdS^{2,1}$ (satisfying the same conditions as in Proposition \ref{principal-curvatures-are-bounded}). For any $p \in S$, let $n_{S}(p)$ denote the unit normal vector to $S$ at $p$, pointing to the future. Let $\nabla$ denote the Levi-Civita connection of $\AdS^{2,1}$. We define the shape operator of $S$ by $B_{S} := -\nabla_{\cdot} n_{S}$, and the third fundamental form by $\III_{S}(\cdot,\cdot) := I_{S}(B_{S}\cdot, B_{S}\cdot).$ (where $I_{S}$ is the induced metric on $S$, which is also called the first fundamental form).\\
We also recall that the product of the eigenvalues of $B_{S}$ (that is $\det(B_{S})$) is equal to $-K_S - 1$, where $K_S$ is the Gaussian curvature of $S$ (which may vary over the surface). Half of the sum of the eigenvalues of $B_{S}$ (that is $\frac{\operatorname{tr}(B_{S})}{2}$, where $\operatorname{tr}(B_{S})$ is the trace of $B_{S}$) is called the mean curvature of $S$ (which may also vary over the surface) and is usually denoted by $H_S$.

\begin{proposition}\label{principal-curvatures-are-bounded}
Let $k > 1$ and $\epsilon > 0$. Let $S \subset \AdS^{2,1}$ be a properly embedded convex space-like disk spanning a $k$-quasicircle $C$. Assume that $S$ is isometric to $(\mathbb{D}, h)$ where $h$ is a conformal complete metric that has curvature in $\left(-\frac{1}{\epsilon}, -1-\epsilon\right)$. Assume that there is a sequence $(M_{p})_{p \in \mathbb{N}}$ of positive numbers such that each derivative of $h$ of order $p$ is bounded by $M_{p}$. Then there exists $D$ that depends on $k$, $\epsilon$, and $(M_p)_{p \in \mathbb{N}}$ such that the principal curvatures of $S$ are in the interval $\left(\frac{1}{D}, D\right)$. 
\end{proposition}

\begin{proof}
We argue by contradiction. Assume the existence of a sequence of properly embedded convex space-like disks $S_{n}$ spanning $k$-quasicircles $C_{n}$, satisfying the hypothesis of Proposition \ref{principal-curvatures-are-bounded}, such that there exists a sequence of points $p_{n} \in S_{n}$ in which one of the principal curvatures at $p_{n}$ goes to infinity. Note that the product of the principal curvatures is bounded (because the curvatures belong to $(-\frac{1}{\epsilon}, -1-\epsilon)$). If one principal curvature goes to $0$, then the other goes to $\infty$. So assume without loss of generality that the largest principal curvature at $p_{n}$ goes to $\infty$.\\
Recall that, by Lemma \ref{compactness-of-quasi-circle}, $C_n$ converge (up to extracting a subsequence) in the Hausdorff topology to a quasi-circle or to the union of a left ruling line and a right ruling line. Up to normalization by isometries, we can assume that $p_{n}$ are equal to a fixed point $p$, and $T_{p_{n}}S_{n}$ are equal to a fixed space-like tangent plane. By Proposition \ref{Hausdorff-converge} and since all $S_n$ are tangent to the same space-like tangent plan, the surfaces $S_{n}$ converge in the Hausdorff sense to a locally convex properly embedded surface $S_{\infty}$ spanning a $k$-quasi-circle $C$.\\
Let $V_{n} : (\mathbb{D},h_n) \to S_{n}$ be isometric embeddings. Let $x_n$ be a sequence of points in the disc $\mathbb{D}$ such that $V_n(x_n) = p$. Let $x_{0} \in \mathbb{D}$ be a fixed point, and let $g_{n}$ be a sequence of elements of $\PSL(2,\mathbb{R})$ such that $g_{n}(x_n) = x_{0}$.\\
We denote $h'_{n} := g_{n}^{*}(h_n)$. Note that the metrics $h'_{n}$ have uniformly bounded derivatives. That is the bounds depend only on the order of the derivatives and not on $n$ or the points, because they are pullbacks of the metrics $h_n$ by hyperbolic isometries (where $h_n$ has uniformly bounded derivatives on the disc at any order).\\
Then, up to extracting a subsequence, $h'_{n}$ converges smoothly ($C^{\infty}$ on compact subsets) to a metric $h'_{\infty}$. By Theorem \ref{JM-thm5.6} and Theorem \ref{JM-nothm5.6}, either $V'_{n} = V_n \circ g_n: (\mathbb{D},h'_{n}) \to S_n$ converge smoothly to an isometric embedding $V'_{\infty}: (\mathbb{D}, h'_{\infty}) \to S_{\infty}$, or there is a complete geodesic $\gamma$ in $(\mathbb{D},h'_{\infty})$ such that $V'_{n}(\gamma)$ converge to a space-like geodesic segment $\Gamma \in \AdS^{2,1}$.\\
We argue by contradiction that we are not in the case of Theorem \ref{JM-nothm5.6}. Suppose this is the case. Then by Theorem \ref{JM-nothm5.6} there is a complete geodesic $\gamma$ in $(\mathbb{D}, h'_{\infty})$ such that $V'_n(\gamma)$ converge to a space-like geodesic $\Gamma$. According to \cite[Lemma 5.4]{schlenker1996surfaces} (and also \cite[Proposition 5.2]{schlenker1996surfaces}) the integral of the mean curvatures of $S_n$ diverges to $+\infty$ at every point in a neighborhood of $\gamma$. Since the products of the eigenvalues of $B_{S_n}$ (where $B_{S_n}$ is the shape operator of $S_n$) are bounded from above and below, the biggest eigenvalue of $B_{S_n}$ diverge to $\infty$ and the smallest eigenvalue converge to $0$. Recall that from \cite[Lemma 5.4]{schlenker1996surfaces} $\gamma$ is limit of curves going in the direction of the smallest eigenvalue. In particular, the length with respect to the third fundamental form of $S_n$ along any geodesic segment transverse to $\gamma$ tends to infinity. By the remark made in \cite{schlenker1996surfaces}, in the paragraphs before and after \cite[Theorem 5.6]{schlenker1996surfaces}, the surfaces $S_n$ accumulate on a pleated surface contained in the past light-like planes directed by the geodesic segment $\Gamma$.\\
On the other hand, by Proposition \ref{Hausdorff-converge}, and since the surfaces $S_n$ have intrinsic curvatures that are uniformly positive and uniformly bounded, it follows that if $C_n$ converges to a $k$-quasi-circle $C$, then the limit surface $S_\infty$ of the sequence $S_n$ cannot be on the past light-like planes directed by the geodesic segment $\Gamma$. This implies that, if we were not in the case of Theorem \ref{JM-thm5.6}, the quasi-circles $C_n$ would have to converge to the union of a left ruling and a right ruling line.\\
Moreover, since all the surfaces $S_n$ are tangent to the same space-like plane at the point $p$, Proposition \ref{Hausdorff-converge} again ensures that their ideal boundaries $C_n$ converge to a $k$-quasi-circle $C$. This contradicts the assumption that the integral of the mean curvatures of $V'_n$ diverges to $+\infty$ at every point in a neighborhood of $\gamma$.\\
We conclude that we cannot be in the case of Theorem \ref{JM-nothm5.6}.\\
Therefore, up to extracting a subsequence, the isometries $V'_n : (\mathbb{D}, h'_n) \to S_n$ converge to an isometric embedding $V'_\infty : (\mathbb{D}, h'_\infty) \to S_\infty$. This shows that the principal curvatures of $S_n$ at $p$ cannot diverge.
\end{proof}
From that, we will deduce that the gluing maps are well defined. Before that, we need to show the following lemma.

\begin{lemma}\label{projection-composed-with-isometries-is-quasi-isometry}
Let $\Omega$ be a globally hyperbolic convex subset spanning a $k$-quasi-circle at infinity. 
Assume that the induced metrics on $\partial^{\pm}\Omega$ have curvatures in 
$\left(-\frac{1}{\epsilon}, -1 - \epsilon\right)$.\\ 
Let $V^{\pm} : (\mathbb{D}, h^{\pm}) \to \partial^{\pm}\Omega$ be isometries. 
Assume that each derivative of $h^{\pm}$ of order $p$ is bounded by $M_p$ on $\mathbb{D}$.\\
Let $\Pi^{\pm}_{l}$ and $\Pi^{\pm}_{r}$ be the left and right projections. Then there exists 
$A > 0$, depending only on $\epsilon$, $k$, and $(M_p)_{p \in \mathbb{N}}$, such that 
$\Pi^{\pm}_{l} \circ V^{\pm}$ and $\Pi^{\pm}_{r} \circ V^{\pm}$ are $A$-quasi-isometries.
\end{lemma}

\begin{proof}
From Proposition \ref{principal-curvatures-are-bounded}, there exists $D > 0$ depending only on $k$, the quasi-symmetric constant of $\partial_{\infty} \Omega$, $\epsilon$, and $(M_p)_{p \in \mathbb{N}}$ such that the principal curvatures of $\partial^{\pm}\Omega$ are in the interval $\left[\frac{1}{D}, D\right]$. By Theorem \ref{estimation-krasnov-schlenker}, we obtain that the projection maps $\Pi_{l}^{\pm}$, $\Pi_{r}^{\pm}$ are $A$-bilipschitz for some constant $A$ that depends on $D$ (which in turn depends on $\epsilon$, $k$, and $(M_p)_{p \in \mathbb{N}}$).\\
Therefore, the maps $\Pi_{l}^{\pm} \circ V^{\pm}: (\mathbb{D},h^{\pm}) \to \mathbb{H}^2$ are $A$-bilipschitz diffeomorphisms. The same argument applies to $\Pi_{r}^{\pm}$.
\end{proof}
From now until the end of the paper, we will always assume that if $\Omega$ is a globally hyperbolic convex subset, then the induced metrics on its boundary components $\partial^{\pm}\Omega$ are isometric to $(\mathbb{D},h^{\pm})$, where $h^{\pm}$ have bounded derivatives of any order on the disc $\mathbb{D}$.\\

Note that from the proof of Lemma \ref{projection-composed-with-isometries-is-quasi-isometry}, the maps $\Pi_{l}^{\pm} \circ V^{\pm}$ extend to a homeomorphism :
$$\partial (\Pi_{l}^{\pm} \circ V^{\pm}): \partial \mathbb{D} \to \partial_{\infty} \mathbb{H}^2.$$

We define $\partial V^{\pm} := (\pi^{\pm}_{l})^{-1} \circ \partial (\Pi_{l}^{\pm} \circ V^{\pm})$, where $\pi^{\pm}_{l}$ is the extension of $\Pi^{\pm}_{l}$ to the boundary at infinity.

\begin{definition}\label{gluing-maps}
Let $\Omega$ be a globally hyperbolic convex subset. Let $V^{\pm}: (\mathbb{D},h^{\pm}) \to \partial^{\pm}\Omega$ be isometries.\\
We define the gluing map to be $\Phi_{\Omega} = (\partial V^{-})^{-1} \circ \partial V^{+}$.\\

\end{definition}
\begin{remark}
Note that the isometries $V^{\pm}$ and the gluing map $\Phi_{\Omega}$ are defined up to composition by M\"obius maps. Also note that M\"obius maps do not necessarily preserve the metrics $h^{\pm}$, so the metrics $h^{\pm}$ are defined up to M\"obius maps. However, we still have a uniquely defined normalized gluing map.\\
Furthermore, if $\partial_{\infty}\Omega$ passes through $0, 1, \infty$, there exist unique metrics $h^{\pm}$ and isometries $V^{\pm}: (\mathbb{D},h^{\pm}) \to \partial^{\pm}\Omega$ such that $\partial V^{\pm}(p) = (p,p)$ for any $p = 0, 1, \infty$.\\
Indeed, if $V^{+}: (\mathbb{D}, h^{+}) \to \partial^{+} \Omega$ and $V'^{+}: (\mathbb{D}, h'^{+}) \to \partial^{+} \Omega$ are two such isometries, then $(V'^{+})^{-1} \circ V^{+}: (\mathbb{D}, h^{+}) \to (\mathbb{D}, h'^{+})$ is an isometry. In particular, since $h'^{+}$ and $h^{+}$ are conformal to $\mathbb{D}$, it follows that $(V'^{+})^{-1} \circ V^{+}$ is also conformal, that is, it belongs to $\PSL(2, \mathbb{R})$. Moreover, since by hypothesis $\partial V^{+}(p) = (p, p)$ and $\partial V'^{+}(p) = (p, p)$ for any $p = 0, 1, \infty$, it follows that $(V'^{+})^{-1} \circ V^{+}$ is the identity because it fixes three points on the boundary of $\mathbb{D}$.    
\end{remark}
\begin{proposition}
 Let $\Omega$ be a globally hyperbolic convex subset. The gluing map is quasi-symmetric.   
\end{proposition}
\begin{proof}
 It follows from the fact that $\Phi_{\Omega} = \partial ((\Pi^{-}_{l} \circ V^{-})^{-1} \circ (\Pi^{+}_{l} \circ V^{+}))$, and each of $\Pi^{-}_{l} \circ V^{-}$ and $\Pi^{+}_{l} \circ V^{+}$ is a bilipchitz diffeomorphism (then quasi-isometric).   
\end{proof}

The discussion above leads to the following proposition.
\begin{proposition}\label{the-siometries-extend-to-the-ideal-boundary}
Let $k > 1$ and $\epsilon > 0$. Let $\Omega$ be a globally hyperbolic convex subset. Assume that $\Omega$ spans a $k$ quasi-circle and assume the existence of the isometries $V^{\pm}: (\mathbb{D},h^{\pm}) \to \partial^{\pm}\Omega$ such that the curvatures of $h^{\pm}$ are in the interval $(-\frac{1}{\epsilon}, -1 - \epsilon)$. Then $V^{\pm}$ extends to a homeomorphism $\overline{V}^{\pm}: \mathbb{RP}^1 \cup \mathbb{D} \to \partial^{\pm}\Omega \cup \partial_{\infty}\Omega$.
\end{proposition}
We denote the gluing map of the convex $\Omega$ by $\Phi_{\Omega}$.\\
We will need the following lemma later.
\begin{lemma}\label{Quasi-isometry-send-to-compact}\cite[Lemma 4.9]{bonsante2021induced}
For any constant $A > 1$ and for any $x \in \mathbb{H}^2$, there exists a compact region $Q$ of $\mathbb{H}^2$ such that if $f$ is a normalized $A$-quasi-isometry of $\mathbb{H}^2$, then $f(x) \in Q$.
\end{lemma}
Now we proceed the proof that the correspondence between the gluing maps and globally hyperbolic convex subsets is proper and continuous.
\begin{proposition}\label{main}
Assume that all the metrics on the disc in the statement are complete and conformal, and they have curvatures in an interval of the form $(-\frac{1}{\epsilon}, -1 - \epsilon)$, for some $\epsilon > 0$. Let $\Omega_{n}$ be a sequence of globally hyperbolic convex subsets such that:
\begin{itemize}
    \item There are isometries $V_{n}^{\pm}: (\mathbb{D}, h_{n}^{\pm}) \to \partial^{\pm}\Omega_{n}$.
    \item The gluing maps $\Phi_{\Omega_{n}}$ are normalized and uniformly quasi-symmetric and they converge in the $C^{0}$ topology to a $k$ quasi-symmetric map $f$.
    \item The metrics $h_{n}^{\pm}$ converge $C^{\infty}$ uniformly on compact subsets to some metrics $h^{\pm}$ on $\mathbb{D}$.
    \item Any derivative of $h^{\pm}_{n}$ or $h^{\pm}$ of order $p$ is bounded by $M_p > 0$ on the disc $\mathbb{D}$ ($M_p$ does not depend on $n$).
\end{itemize}
Then globally hyperbolic convex subsets $\Omega_{n}$ converge to a globally hyperbolic convex subset $\Omega$ in the Hausdorff sense, and the isometries $V_{n}^{\pm}: (\mathbb{D}, h_{n}^{\pm}) \to \partial^{\pm}\Omega_{n}$ converge to isometries $V^{\pm}: (\mathbb{D}, h^{\pm}) \to \partial^{\pm}\Omega$, and $\Phi_{\Omega}$, the gluing map of $\Omega$, is equal to $f$.
\end{proposition}
We will split the proof of Proposition \ref{main} into lemmas and propositions.\\
\subsection{Properness of the map $\Phi_{\cdot}$ }
In this section we prove the properness of the  map $\Phi_{\cdot}$ in the sense of the following proposition.
\begin{proposition}\label{proper}
Let $\epsilon > 0$ and $(M_{p})_{p \in \mathbb{N}}$ a sequence of positive numbers. For any $k > 1$ there exists $k' > 1$ such that for any globally hyperbolic convex subset $\Omega$ in which the induced metrics on $\partial^{\pm}\Omega$ have curvatures in $(-\frac{1}{\epsilon},-1-\epsilon)$ and any derivative of the metric of order $p$ is bounded by $M_p$. We assume that the gluing map is normalized. If the gluing map $\Phi_{\Omega}$ is $k$-quasi-symmetric then $\partial_{\infty} \Omega$ is a $k'$ quasi-circle.   
\end{proposition}
Before proceeding with the proof, let us give the following lemma.
\begin{lemma}\cite[Lemma 8.4]{bonsante2021induced}\label{planes-bdms}
Let $P$ be a totally geodesic space-like plane. We denote by $P^{+}$ (resp $P^{-}$) the union of all future-oriented (resp past-oriented ) time-like geodesic segments of length $\frac{\pi}{2}$ starting orthogonally from $P$.  
\begin{itemize}
    \item let $S \subset \AdS^{2,1}$ be a space-like past convex surface, and let $P$ be a space-like totally geodesic plane. In the neighborhood of all points $x \in S \cap P$ , the intersection $S \cap P^{+}$ is locally convex in the induced metric on $S$.
    \item Let $S \subset \AdS^{2,1}$ be a space-like future convex surface, and let $P$ be a space-like totally geodesic plane. In the neighborhood of all points $x \in S \cap P$ , the intersection $S \cap P^{-}$ is locally convex in the induced metric on $S$.
\end{itemize}    
\end{lemma}
Now we give a proof of Proposition \ref{proper}, note that the proof is similar to the proof of \cite[Proposition 8.3]{bonsante2021induced}.
\begin{proof}
We argue by contradiction. Suppose that $\Omega_n$ is a sequence of normalized globally hyperbolic convex subsets such that $k'_{n}$, the optimal quasi-symmetric constant of the ideal boundary $\partial_{\infty} \Omega_n$, diverges to infinity. We will prove that $k_{n}$, the quasi-symmetric constants of their gluing maps $\Phi_{\Omega_{n}}$, must also go to $\infty$. Let $V^{\pm}_{n}:(\mathbb{D},h^{\pm}_{n}) \to \partial^{\pm}\Omega_{n}$ be isometries such that the normalized gluing map $\Phi_{\Omega_{n}} = (\partial V^{-}_{n})^{-1} \circ (\partial V^{+}_{n})$. By applying an isometry of $\AdS^{2,1}$ to $\Omega_{n}$, we can assume that $\partial V^{+}_{n}(p) = (p,p)$ for any $p = 0,1,\infty$. Since we assumed that the gluing map is normalized, it follows that $\partial V^{-}_{n}(p) = (p,p)$ for any $p = 0,1,\infty$. Given $k'_{n} \to \infty$, Proposition \ref{rhombus} implies that the width $w(C_{n})$ of $C_{n} := \partial_{\infty}\Omega_{n}$ goes to $\frac{\pi}{2}$. After adjusting by isometries, we can assume that $C_{n}$ converges to a rhombus $C_{\diamond}$ as in \cite[Example 6.7]{bonsante2021induced} (note that even after applying these isometries, we can still assume that $C_{n}$ and $C_{\diamond}$ are normalized). Recall that under this normalization, the width of $C_{\diamond}$ is $\frac{\pi}{2}$.\\
Let's work in the projective model of $\AdS^{2,1}$ with the coordinates so that in the affine chart $x_{4} = 1$, $\partial_{\infty}\AdS^{2,1}$ is the hyperboloid $x_{1}^{2} + x_{2}^{2} = x_{3}^{2} + 1$, with $\AdS^{2,1}$ seen as the region $x_{1}^{2} + x_{2}^{2} < x_{3}^{2} + 1$. We may then arrange that the vertices of $C_{\diamond}$ are the points $(\pm \sqrt{2},0,-1)$ and $(0,\pm\sqrt{2},1)$. Since $\partial^{+}\Omega_{n}$ is in the future of $\CH(C_{n})$ but contained in $D(C_n)$, the domain of dependence of $C_n$, and since both $\overline{D(C_{n})}$ (the closure of $D(C_{n})$ in $\AdS^{2,1} \cup \partial_{\infty}\AdS^{2,1}$) and $\CH(C_{n})$ converge to $\CH(C_{\diamond})$, we have that $\partial^{+} \Omega_{n}$ converge to $S^{+}_{\infty}$, where $S^{+}_{\infty}$ is the union of the two future faces of $\CH(C_{\diamond})$. Similarly, $\partial^{-} \Omega_{n}$ converge to $S^{-}_{\infty}$, where $S^{-}_{\infty}$ is the union of the two past faces of $\CH(C_{\diamond})$.\\
Let $P$ be the space-like totally geodesic plane equal to the intersection with $\AdS^{2,1}$ of the plane $x_3 = 0$ in $\mathbb{R}^3$. Let $a = \left(-\frac{\sqrt{2}}{2}, \frac{\sqrt{2}}{2}, 0 \right)$, $b = \left(-\frac{\sqrt{2}}{2}, -\frac{\sqrt{2}}{2}, 0 \right)$, $c = \left(\frac{\sqrt{2}}{2}, -\frac{\sqrt{2}}{2}, 0 \right)$, and $d = \left(\frac{\sqrt{2}}{2}, \frac{\sqrt{2}}{2}, 0 \right)$. Thus, $a, b, c, d$ are the intersection points of $\partial P$ with $C_{\diamond}$ occurring in the cyclic order $a, b, c, d$.\\ 
For all $n \in \mathbb{N}$, choose intersection points $a_{n}, b_{n}, c_{n}, d_{n}$ of $C_{n}$ with $P$ such that $a_{n} \to a$, $b_{n} \to b$, $c_{n} \to c$, and $d_{n} \to d$. Such points exist since $C_n$ converge to $C_{\diamond}$. Define $a_{n}^{\pm}, b_{n}^{\pm}, c_{n}^{\pm}, d_{n}^{\pm}$ by the equalities $\partial V^{\pm}_{n}(a^{\pm}_{n}) = a_{n}$, $\partial V^{\pm}_{n}(b^{\pm}_{n}) = b_{n}$, $\partial V^{\pm}_{n}(c^{\pm}_{n}) = c_{n}$, $\partial V^{\pm}_{n}(d^{\pm}_{n}) = d_{n}$. Note that $a^{-}_{n}, b^{-}_{n}, c^{-}_{n}, d^{-}_{n}$ are the images of $a^{+}_{n}, b^{+}_{n}, c^{+}_{n}, d^{+}_{n}$ under $\Phi_{\Omega_{n}}$.\\ 
We will show that the cross-ratio of $a^{+}_{n}, b^{+}_{n}, c^{+}_{n}, d^{+}_{n}$ tends to $0$, while the cross-ratio of $a^{-}_{n}, b^{-}_{n}, c^{-}_{n}, d^{-}_{n}$ tends to infinity.\\
Let $P^{+}$ be the future of $P$ as in Lemma \ref{planes-bdms}. Let $Q$ be the time-like plane defined by the equation $x_{2} = 0$. Then the path $Q \cap S^{+}_{\infty} \cap P^{+}$ from $Q \cap \overline{ab}$ to $Q \cap \overline{cd}$ along the piecewise light-like geodesic $Q \cap S^{+}_{\infty}$ has length zero in the $\AdS^{2,1}$ metric (here, $\overline{ab}$ and $\overline{cd}$ are the geodesics inside $S_{\infty}^{+}$). This implies that the lengths of the paths $Q \cap S_{n}^{+} \cap P^{+}$ converge to zero.\\
Note that $S^{+}_{n} \cap P^{+}$ has a locally convex boundary with respect to the induced metric of $S_{n}^{+}$. Then the set $U^{+}_{n} := (V^{+}_{n})^{-1}(S^{+}_{n} \cap P^{+})$ is a region of $(\mathbb{D}, h_{n}^{+})$ that has a locally convex boundary, so it is globally convex (because the metrics are non positively curved). Since $a_{n}, b_{n}, c_{n}, d_{n}$ are the intersection points of $C_{n}$ and $P^{+}$ by Proposition \ref{the-siometries-extend-to-the-ideal-boundary}, the set $U^{+}_{n}$ contains the points $(\partial V_{n}^{+})^{-1}(a_{n}) = a_{n}^{+}$, $(\partial V_{n}^{+})^{-1}(b_{n}) = b_{n}^{+}$, $(\partial V_{n}^{+}(c_{n}))^{-1} = c_{n}^{+}$, $(\partial V_{n}^{+}(d_{n}))^{-1} = d_{n}^{+}$ in its ideal boundary. Also, since $U_{n}^{+}$ is convex, it contains the geodesics (with respect to the metric $h_{n}^{+}$) $\overline{a_{n}^{+}b_{n}^{+}}$ and $\overline{c_{n}^{+}d_{n}^{+}}$, which implies that $S_{n}^{+} \cap P^{+}$ contains the geodesics $\gamma(\overline{a_{n}^{+}b_{n}^{+}}) := V_{n}^{+}(\overline{a_{n}^{+}b_{n}^{+}})$ and $\gamma(\overline{c_{n}^{+}d_{n}^{+}}) := V_{n}^{+}(\overline{c_{n}^{+}d_{n}^{+}})$ (in other words, $\gamma(\overline{a_{n}^{+}b_{n}^{+}})$ and $\gamma(\overline{c_{n}^{+}d_{n}^{+}})$ are the geodesics in $S_{n}^{+}$ that are the images of $\overline{a_{n}^{+}b_{n}^{+}}$ and $\overline{c_{n}^{+}d_{n}^{+}}$, respectively, under $V_{n}^{+}$).\\
Also note that the path $Q \cap S_{n}^{+} \cap P_{+}$ crosses the two geodesics $\gamma(\overline{a_{n}^{+}b_{n}^{+}})$ and $\gamma(\overline{c_{n}^{+}d_{n}^{+}})$. Recall that the length of $Q \cap S_{n}^{+} \cap P_{+}$ converges to $0$, so the distance between the two geodesics $\gamma(\overline{a_{n}^{+}b_{n}^{+}})$ and $\gamma(\overline{c_{n}^{+}d_{n}^{+}})$ goes to zero. It follows that the cross-ratio of $(a^{+}_{n}, b^{+}_{n}, c^{+}_{n}, d^{+}_{n})$ goes to $0$.\\
We apply a similar argument on the surfaces $S_{n}^{-}$. We denote by $T$ the time-like plane defined by the equation $x_1 = 0$. The path $T \cap S_{\infty}^{-} \cap P^{-}$ has length zero. By similar arguments as above (applying Lemma \ref{planes-bdms} on $P^{-} \cap S_{\infty}^{-}$ and considering the path $T \cap S^{-}_{\infty} \cap P^{-}$), the cross-ratio $(d_{n}^{-},a_{n}^{-},b_{n}^{-},c_{n}^{-})$ converges to $0$. This implies that the cross-ratio of $(a_{n}^{-},b_{n}^{-},c_{n}^{-},d_{n}^{-})$ converges to $\infty$. We deduce that the quasi-symmetric constants of the gluing maps $\Phi_{\Omega_{n}}$ diverge.
  
\end{proof}
\subsection{Continuity of the map $\Phi_{\cdot}$}

The proof that we give for the next proposition is similar to the proof of \cite[Proposition 8.2]{bonsante2021induced}.
\begin{proposition}\label{continuty-of-the gluing-maps}
Let $\Omega_{n}$ be a sequence of globally hyperbolic convex subsets. Assume the existence of isometries $V_{n}^{\pm}:(\mathbb{D},h_{n}^{\pm}) \to \partial^{\pm}\Omega_{n}$ such that the metrics $h^{\pm}_{n}$ converge $C^{\infty}$ on compact subsets to $h^{\pm}$, and all the metrics have curvatures in $(-\frac{1}{\epsilon}, -1-\epsilon)$ for some $\epsilon > 0$. Assume there exists a sequence of positive numbers $(M_p)_{p \in \mathbb{N}}$ such that any derivative of $h^{\pm}_n$ or $h^{\pm}$ of order $p$ is bounded by $M_p$. Also assume that all the gluing maps $\Phi_{\Omega_{n}}$ are normalized and uniformly quasi-symmetric (they are all $k$-quasi-symmetric for the same $k$), and that $\Phi_{\Omega_{n}}$ converge to a quasi-symmetric map $f$ in the $C^{0}$ topology. Assume further that $\Omega_{n}$ converge (which is always possible, up to extracting a subsequence, after applying isometries on $\Omega_{n}$) in the Hausdorff topology to a hyperbolic convex subset $\Omega$.\\
Then $V_{n}^{\pm}:(\mathbb{D},h_{n}^{\pm}) \to \partial^{\pm}\Omega_{n}$ converge to isometries $V^{\pm}:(\mathbb{D},h^{\pm}) \to \partial^{\pm}\Omega$, and the gluing map $\Phi_{\Omega}$ of $\Omega$ is equal to $f$.
\end{proposition}
\begin{proof}
Up to normalizing by isometries of $\AdS^{2,1}$, we can assume that $\partial V^{+}_{n}(p) = (p,p)$ for any $p = 0,1,\infty$. Note that after this normalization, $\Omega_{n}$ converges in the Hausdorff topology to some convex subset $\Omega$ that has a quasi-circle as ideal boundary. Since we assumed that the gluing maps $\Phi_{\Omega_{n}}$ are normalized, we also have that $\partial V^{-}_{n}(p) = (p,p)$ for any $p = 0,1,\infty$.\\
We will show that each of the isometries $V_{n}^{\pm}: (\mathbb{D}, h_{n}^{\pm}) \to \partial^{\pm} \Omega_{n}$ converges to an isometry $V^{\pm}: (\mathbb{D}, h^{\pm}) \to \partial^{\pm} \Omega$. To do that, we will show that there exists $x_{0} \in \mathbb{D}$ such that $V^{\pm}_{n}(x_{0})$ lie in a compact subset of $\AdS^{2,1}$. Also recall that, by Proposition \ref{principal-curvatures-are-bounded}, the principal curvatures of $\partial^{\pm} \Omega_n$ are uniformly bounded. So the convergence of the isometries $V_{n}^{\pm}$ will follow from Theorem \ref{JM-thm5.6}.\\
We denote by $\Pi^{\pm}_{n,l}$ (resp $\Pi^{\pm}_{n,r}$) the left (resp the right) projection from $\partial^{\pm}\Omega_n$. Let $x_{0} \in \mathbb{D}$ be a fixed point. From Lemma \ref{projection-composed-with-isometries-is-quasi-isometry}, there is $A$ that depends only of $\epsilon$, $k$ and the sequence $(M_p)_{p \in \mathbb{N}}$ such that each of $\Pi^{\pm}_{n,l} \circ V^{\pm}_{n}$ and $\Pi^{\pm}_{n,r} \circ V^{\pm}_{n}$ is an $A$ quasi-isometry from $\mathbb{H}^{2}$ to $\mathbb{H}^2$. Also, since $\partial V^{\pm}_{n}(p) = (p,p)$, for any $p = 0,1,\infty$, the $A$ quasi-isometries $\Pi^{\pm}_{n,l} \circ V^{\pm}_{n}$ and $\Pi^{\pm}_{n,r} \circ V^{\pm}_{n}$ are normalized. Then by Lemma \ref{Quasi-isometry-send-to-compact} there are compact subsets $K_{l}$ and $K_{r}$ of $\mathbb{H}^2$ such that the images $\Pi^{\pm}_{n,l} \circ V^{\pm}_{n}(x_{0})$ and $\Pi^{\pm}_{n,r} \circ V^{\pm}_{n}(x_{0})$ belong to $K_{l}$ and $K_{r}$ respectively.\\
Hence $V^{\pm}_{n}(x_{0})$ lie on a subset of time-like geodesics $L_{x,y}$ where $x$ vary in $K_{l}$, a compact subset of $\mathbb{H}^2$, and $y$ vary in $K_{r}$, also  a compact subset of $\mathbb{H}^2$. This implies that each of the sequences $V^{\pm}_{n}(x_{0})$ lies on a compact subset of $\AdS^{2,1}$ (recall that $\Omega_{n}$ converge in the Hausdorff sense). Then the isometries $V^{\pm}_{n}: (\mathbb{D},h^{\pm}_{n}) \to \partial^{\pm}\Omega_{n}$ converge to the isometry $V^{\pm}: (\mathbb{D},h^{\pm}) \to \partial^{\pm}\Omega$. Note that the gluing map $\Phi_{\Omega_{n}}$ satisfies:
$$\Phi_{\Omega_{n}} = \partial ((\Pi^{-}_{n,l} \circ V^{-}_{n})^{-1} \circ (\Pi^{+}_{n,r} \circ V^{+}_{n})),$$ 
and note that all maps $(\Pi^{-}_{n,l} \circ V^{-}_{n})^{-1} \circ (\Pi^{+}_{n,r} \circ V^{+}_{n})$ are normalized uniformly quasi-isometries. Since $(\Pi^{-}_{n,l} \circ V^{-}_{n})^{-1} \circ (\Pi^{+}_{n,r} \circ V^{+}_{n})$ converges to $(\Pi^{-}_{l} \circ V^{-})^{-1} \circ (\Pi^{+}_{r} \circ V^{+})$, then also $\Phi_{\Omega_{n}}$ converge to $\partial ((\Pi^{-}_{l} \circ V^{-})^{-1} \circ (\Pi^{+}_{r} \circ V^{+})) = (\partial V^{-})^{-1} \circ \partial V^{+} = \Phi_{\Omega} $. Since we have also assumed that $\Phi_{\Omega_{n}}$ converges to $f$, we deduce $\Phi_{\Omega} = f$.

\end{proof}
\section{Approximation}
In the proof of our main theorem, we will use an approximation by lifts of globally hyperbolic manifolds. For that purpose, we need to approximate the given metrics on the disc by metrics that are invariant under the action of Fuchsian groups.
\begin{proposition}\label{approx-the-metrics}
Let $S_{n}$ be a sequence of closed surfaces having genus $g_{n}$ converging to $\infty$. Let $\rho_{n}: \pi_{1}(S_{n}) \to \PSL(2,\mathbb{R})$ be a sequence of Fuchsian representations that have injectivity radius going to $\infty$. Let $h$ be a complete, conformal Riemannian metric on the disc $\mathbb{D}$ that has curvature in $\left(-\frac{1}{\epsilon}, -1-\epsilon\right)$, for some $\epsilon > 0$. also assume that there is a sequence of positive real numbers $(M_p)_{p \in \mathbb{N}}$ such that each derivative of $h$ of order $p$ is bounded by $M_p$. Then there exists a sequence $h_{n}$ of complete, conformal Riemannian metrics that have curvature in $\left(-\frac{1}{\epsilon'},-1-\epsilon'\right)$, for some $\epsilon' > 0$, each $h_{n}$ is $\rho_{n}$-invariant, and $h_{n}$ converge $C^{\infty}$ uniformly on compact subsets of $\mathbb{D}$ to $h$. Moreover, there is a sequence of positive real numbers $(M'_p)_{p \in \mathbb{N}}$ such that each derivative of $h_{n}$ of order $p$ is bounded by $M'_p$.
\end{proposition}
To do that, we will use curvatures. We will approximate $K_{h}$, the curvature of $h$, by smooth functions $K_{n}$ invariant under $\rho_{n}$. Then, we will show that each $K_{n}$ is the curvature of a metric $h_{n}$ which is invariant under $\rho_{n}$. Finally, we will show that the metrics $h_{n}$ converge to $h$ smoothly on compact subsets of $\mathbb{D}$. Note that for us it is important to have bounded derivatives with bounds that do not depend on $n$.
\begin{lemma}\label{construct-curvatures}
Let $K : \mathbb{D} \to \left(-\frac{1}{\epsilon},-1-\epsilon\right)$ be a smooth function such that any derivative of it at order $p$ is bounded by some $M_p > 0$ uniformly on $\mathbb{D}$. Let $\rho_{n} : \pi_{1}(S_{n}) \to \PSL(2,\mathbb{R})$ be a sequence of Fuchsian representations that have injectivity radius growing to $\infty$. Then there exists a sequence of smooth functions $K_{n} : \mathbb{D} \to \left(-\frac{1}{\epsilon'},-1-\epsilon' \right)$, such that each $K_{n}$ is $\rho_{n}$-equivariant, $K_{n}$ converge $C^{\infty}$ on compact subsets to $K$, and each derivative of order $p$ of $K_{n}$ is bounded on the disc $\mathbb{D}$ by some $M'_p$, where $M'_p$ does not depend on $n$.
\end{lemma}
\begin{proof}
Let $D_{n}$ be a fundamental domain of $\rho_{n}$. Let $0 < r_{n}$ be such that $B(o,2r_{n}) \subset D_{n}$, where $B(o,r)$ is the hyperbolic ball centered at $o$, the center of $\mathbb{D}$. Since the injectivity radius of $\rho_{n}$ go to $\infty$, we can assume that $r_{n}$ go to $\infty$.\\
Let $\phi: \mathbb{R} \to  \left [0,1  \right ]$ be a smooth function such that $\phi \mid_{\left (-\infty,1  \right ]} = 1$ and $\phi \mid_{\left [2,+\infty  \right )} = 0$ (see Urysohn smooth lemma). We define the function $K'_n$ on $D_n$ by :  
$$
K'_n(x) := \phi\left(\frac{d_{\mathbb{H}^2}(o,x)}{r_n}\right) K + \left(1 - \phi\left(\frac{d_{\mathbb{H}^2}(o,x)}{r_n}\right)\right) \left(-\frac{1}{\epsilon}\right).
$$
In particular, note that $K'_n\mid_{B(o,r_n)} = K$ and $K'_n\mid_{D_n \setminus B(o,2r_n)} = -\frac{1}{\epsilon}$.\\
Note that that $K'_n$ has values in the interval $\left [-\frac{1}{\epsilon},-1-\epsilon  \right ]$. Also note that the derivatives of $K'_{n}$ depend only on $K$ and its derivatives, the derivatives of $\phi$ on the closed interval $\left[1, 2\right]$, and the derivatives of $\frac{d_{\mathbb{H}^2}(o,x)}{r_n}$ that are uniformly bounded at any order $\alpha$ on a neighborhood of the subset $B(o,2r_n) \setminus B(o,r_n)$. Then we deduce that any derivative of $K'_{n}$ of order $\alpha$ is bounded by a constant $M'_{\alpha}$ that does not depend on $n$ (it depends on the function $\phi$ and on the bound of the derivative of $K$ of order $\alpha$).\\
We define the map $K_{n}$ on the disc $\mathbb{D}$ by extending $K'_{n}$ by reflections (that is $K_{n}$ is the unique $\rho_{n}$ invariant map that extends $K'_{n}$ to the disc). Note that $K_{n}$ is smooth since $K'_n$ is constant on a neighborhood of the boundary of the fundamental domain $D_{n}$. It is also clear that $K_{n}$ converge uniformly on compact subsets to $K$ (because $r_{n}$ is growing to $\infty$). 
\end{proof}

We will use the next theorem : 
\begin{theorem}\label{compact}\cite{kazdan1974curvature}
Let $K: S \to \mathbb{R}_{-}$ be a $C^{\infty}$ function, and let $ \left [ g \right ]$ be a conformal class on $S$. Then there exists a unique complete metric $h$ on $S$ conformal to $g$ that has curvature equal to $K$.    
\end{theorem}
\begin{lemma}\label{nathanial}
Let $h_{n}$ be a sequence of complete metrics on $\mathbb{D}$, and let $h$ be also a complete metric on $\mathbb{D}$. Assume that all the metrics are conformal to $ \left|dz \right|^{2}$.\\
For each $n$ we denote by $K_n$ the curvature of $h_n$. Moreover assume that $K_{h}$ the curvature of $h$, and $K_n$ for any $n$ belong to $\left [ -\frac{1}{\epsilon},-1-\epsilon \right ]$ for some $\epsilon > 0$.\\
If $K_{n}$ converge uniformly $C^{\infty}$ on compact subsets to $K_{h}$, then $h_{n}$ converge uniformly $C^{\infty}$ on compact subsets to $h$.\\
Moreover, if there is a sequence $(M_p)_{p \in \mathbb{N}}$ such that any derivative of $K_{n}$ of order $p$ is bounded by $M_p$, then there is a sequence of positive real numbers $M'_p$ such that any derivative of $h_n$ of order $p$ is bounded by $M'_p$.
\end{lemma}
\begin{proof}
Let $f$ be a smooth function on a closed hyperbolic disc $B$.\\
In this proof we will use euclidean norms, and we will use the fact that the hyperbolic metric and the euclidean metric are bi-Lipchitz on compact subsets.\\
We denote :
$$\left\|f \right\|_{k,p} := (\sum_{\alpha \leq k} \int_{B} \left|\partial^{\alpha}f \right|^{p}dx )^{\frac{1}{p}},$$
also we denote :
$$ \left\|f \right\|_{C^{k,\alpha}(B)} := \sum_{\left|\beta \right| \leq k}  \left\|\partial^{\beta}f \right\|_{\infty} + \sum_{\left|\beta \right| = k} \left [ \partial^{\beta}f \right ]_{\alpha},$$
where $\left [ \partial^{\beta}f \right ]_{\alpha} := \sup\limits_{0 < R < 1, z \in B } R^{-\alpha}\sup\left\{\left|f(x)-f(y) \right|; x,y \in D(z,R) \cap B \right\}$ (here, by $D(z,R)$, we mean the euclidean ball).\\
Recall that at a point $p$ of the disc $\mathbb{D}$, the hyperbolic metric $h_{-1}$ is defined by $h_{-1} := \frac{ 4 \left|dz \right|^{2}}{(1-\left\|p \right\|^{2})^{2}}$.\\
We will show that if $K$, the curvature of a complete conformal metric $h$ on $\mathbb{D}$, belongs to an interval of the form $(-\frac{1}{\epsilon},-1-\epsilon)$ and every derivative of $K$ or order $p$ is bounded by $M_p > 0$, then all the partial derivatives of the metric at any order are bounded by constants that depend only on $M_p$ and $\epsilon$.\\
Fix $r > 0$.  Denote by $B(x,r)$ the hyperbolic ball of center $x$ and radius $r$. The metric $h$ satisfies $h = e^{2u}h_{-1}$, where $u$ is a smooth function and $h_{-1}$ is the hyperbolic metric. Note that (see \cite{troyanov1991schwarz}) the function $u$ is bounded on the disc because the curvature of $h$ is negative and bounded.\\
Up to applying a hyperbolic isometry that sends the ball $B(x, r)$ to the hyperbolic ball $B(o, r)$, we may assume that we remain in the hyperbolic ball $B := B(o, r)$ centered at $o$ with hyperbolic radius $r$. Let $g$ be a hyperbolic isometry, then we get $g^{*}(h) = e^{u^{*}}h_{-1}$, where $u^{*} = u \circ g$, in particular $u^{*}$ still bounded on $\mathbb{D}$.\\
Since the hyperbolic metric and the euclidean metric are bi-Lipchitz with a bi-Lipchitz constant that depends only on $r$, we find that all the derivatives of $K$ restcited to $B$ are bounded with respect to the euclidean metric at any order. We will show that this implies that the derivatives of $u^{*}$ are bounded by constants that depend only on the bounds of the derivatives of $K$. Again, since the euclidean metric and the hyperbolic metric are bi-Lipchitz with a constant that depends only on $r$ (that we can control), then the derivatives of $u^{*}$ will be bounded by the hyperbolic metrics and the bounds depend only on $r$ and the hyperbolic bounds of the derivatives of $K$. After applying an isometry that maps the ball $B$ back to the ball $B(x, r)$, we find that all derivatives of $u$ are bounded on $B(x, r)$ with respect to the hyperbolic metric. These bounds depend only on $\epsilon$ and the bounds of $K$.
\\   
 Now we restrict our self to the hyperbolic ball $B(o,r)$. Recall the equation $e^{2u^{*}}K^{*}+1 = -\Delta u^{*}$, where $h^{*} = e^{2u^{*}}h_1$ and $\Delta u^{*}$ is the hyperbolic Laplacian and $K^{*}$ is the curvature of $h^{*}$. By \cite[Theorem 10.3.1]{nicolaescu2020lectures}, for any $k \in \mathbb{N}$ and $p > 1$ there exists a constant $C_{k,p}$ such that the following inequality holds :
 $$  \left\|u^{*} \right\|_{k+2,p} \leq C_{k,p}( \left\|\Delta u^{*} \right\|_{k,p}+ \left\|u^{*} \right\|_{p}).$$
 Note that since $u^{*}$ is uniformly bounded with bounds that depend only  $\epsilon$, and since $K^{*}$ and all its derivatives are uniformly bounded at any order, it follows that by applying induction on the equality $e^{2u^{*}}K^{*}+1 = -\Delta u^{*}$ we get the existence of constants $M'_{k,p}$ such that  $\left\|u^{*} \right\|_{k,p}$ is bounded for any $k,p$ by $M'_{k,p}$.\\
 By Morrey inequality there exist constant $C'_{k,p}$ such that :
 $$\left\|u^{*} \right\|_{C^{k-\gamma,\gamma-\frac{2}{p}}} \leq C'_{k,p}\left\|u^{*} \right\|_{k,p},$$
 where $\gamma = \left \lfloor \frac{2}{p}\right \rfloor+1$.\\
 Then there exists $M''_{k,\alpha}$ that depend only on $M_p$ and $\epsilon$ such that $\left\|u^{*} \right\|_{C^{k,\alpha}(B(x,r))}$ is bounded by $M''_{k,\alpha}$, this implies in particular that all the derivatives of $u^{*}$ are bounded at any order by constants that depends only on $\epsilon$ and $(M_p)_{p \in \mathbb{N}}$.\\
 In particular, the lemma will follows because the metrics $h_{n}$ will converge uniformly $C^{\infty}$ on compact subsets up to extracting a subsequence (because their derivatives are uniformly bounded at any order with respect the hyperbolic metric, this follows from the arguments above). Since the curvatures of $h_{n}$ converge to the curvature of $h$, we obtain that $h_{n}$ converge (up to extracting) to $h$.\\
\end{proof}
Now we give a proof for Proposition \ref{approx-the-metrics}.
\begin{proof}
By Lemma \ref{construct-curvatures} we construct $K_{n}: \mathbb{D} \to (-\frac{1}{\epsilon'},-1-\epsilon')$, a sequence of $\rho_{n}$ invariant functions that converge $C^{\infty}$ on compact subsets to $K_{h}$, where $K_{h}$ is the curvature of $h$ (recall that $K_{n}$ have uniformly bounded derivatives). By Theorem \ref{compact}, for any $n$ there is a unique complete conformal metric $h_{n}$ on the disc which has curvature equal to $K_{n}$, and which is $\rho_{n}$ invariant. Since $K_{n}$ the curvatures of the metrics $h_{n}$ converge smoothly to $K_{h}$, Lemma \ref{nathanial} implies that the metrics $h_{n}$ converge, up to extracting a subsequence, to $h$. The derivatives of $h_n$ at any order are all uniformly bounded (independently on $n$) by Lemma \ref{nathanial}.
\end{proof}
\section{The proof of the main theorem}
\begin{figure}
    \centering
    \includegraphics[width=4cm]{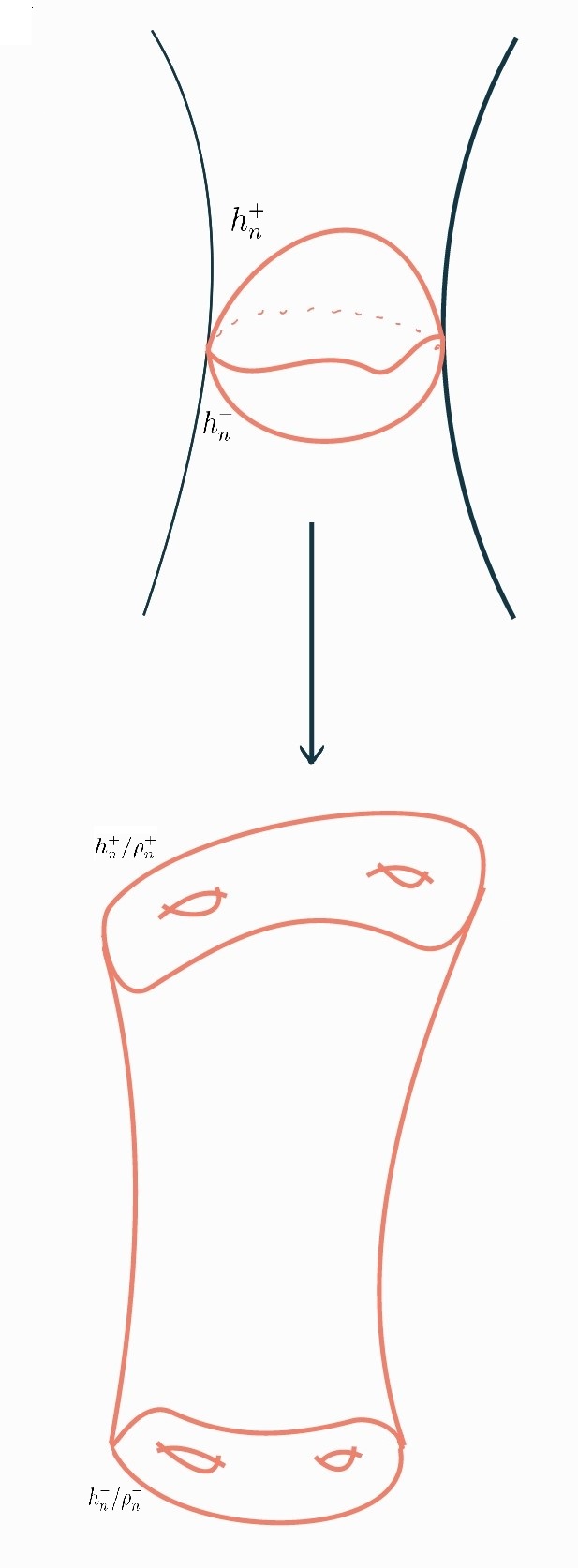}
    \caption{If there is a closed hyperbolic surface $S_n$, two Fuchsian representations $\rho_{n}^{+}$ and $\rho_{n}^{-}: \pi_{1}(S_n) \to \PSL(2,\mathbb{R})$, and two complete Riemannian metrics $h^{+}_{n}$ and $h^{-}_{n}$ on $\mathbb{D}$ that have sectional curvature strictly less than $-1$ and are invariant under the action of $\rho_{n}^{+}$ and $\rho_{n}^{-}$ respectively. Then by Theorem \ref{tumb}, we can find a globally hyperbolic manifold diffeomorphic to $M_{n} = S_n \times [0,1]$ such that the induced metric on $S_n \times  \left\{ 0\right\}$ is homotopic to $h^{+}_{n} / \rho_{n}^{+}$, and the induced metric on $S_n \times  \left\{ 1\right\}$ is homotopic to $h^{-}_{n} / \rho_{n}^{-}$. The globally hyperbolic manifolds $M_{n}$ lift to globally hyperbolic convex subsets $\Omega_{n}$. There is a unique quasi-symmetric map equivariant under the action of $\rho^{+}_{n}$ and $\rho^{-}_{n}$, the gluing map of $\Omega_{n}$ must be equal to this map. The convex sets $\Omega_{n}$ will converge to the desired $\Omega$ that we want to realize.}
    \label{tamburelli-pic}
\end{figure}
In this section we will prove the main theorem (see Figure \ref{tamburelli-pic}). Before proceeding with the proof we need to recall the following statements.\\
The next theorem, is the group action invariant version of our main theorem.

\begin{theorem}\cite{tamburelli2018prescribing}\label{tumb}
Given two metrics $g^+$ and $g^-$ with curvature $\kappa < -1$ on a closed, oriented surface $S$ of genus $g \geq 2$, there exists an $\AdS^{2,1}$ manifold $N$ with smooth, space-like, strictly convex boundary such that the induced metrics on the two connected components of $\partial N$ are isotopic to $g^+$ and $g^-$.
\end{theorem}
The next proposition shows that equivariant quasi-symmetric maps are dense in the set of quasi-symmetric maps.
\begin{proposition}\cite[Proposition 9.1]{bonsante2021induced}\label{approx-the-quasi}
Let $f$ be a normalized quasi-symmetric map. There is a sequence of equivariant normalized uniformly quasi-symmetric maps, $\rho_{n}^{+}, \rho_{n}^{-}: \pi_{1}(S_{n}) \to \PSL(2,\mathbb{R})$, that converge to $f$. Here, $S_{n}$ is a sequence of closed surfaces with genus $g_{n}$ going to $\infty$, and $\rho_{n}^{+}, \rho_{n}^{-}$ are a sequence of Fuchsian representations whose injectivity radius go to $\infty$.
\end{proposition}
Now we proceed the proof of our main theorem.
\begin{theorem}
Let $h^{\pm}$ be two complete conformal Riemannian metrics on the disk $\mathbb{D}$ that have curvature in $(-\frac{1}{\epsilon}, -1-\epsilon)$ for some $\epsilon > 0$, and each derivative of $h^{\pm}$ of order $p$ is bounded by $M_p > 0$. Let $f$ be a normalized quasi-symmetric map. There exists a normalized globally hyperbolic convex subset $\Omega$, and normalized isometries $V^{\pm}: (\mathbb{D}, h^{\pm}) \to \partial^{\pm} \Omega$ such that $f = \Phi_{\Omega}$.
\end{theorem}
\begin{proof}
By Proposition \ref{approx-the-quasi}, there are quasi-Fuchsian representations $\rho_{n}^{+}, \rho_{n}^{-}: \pi_{1}(S_{n}) \to \PSL(2,\mathbb{R})$, and $f_{n}$ a sequence of $\rho_{n}^{+}, \rho_{n}^{-}$ equivariant quasi-symmetric maps that converge in the $C^{0}$ topology to $f$. Also by Proposition \ref{approx-the-metrics} there is a $\rho_{n}^{+}$ (resp $\rho_{n}^{-}$) invariant metrics $h_{n}^{+}$ (resp $h_{n}^{-}$) that converge $C^{\infty}$ on compact subsets of $\mathbb{D}$ to $h^{+}$ (resp $h^{-}$). Also by Proposition \ref{approx-the-metrics}, the metrics $h^{\pm}_{n}$ have uniformly bounded derivatives.\\ 
By Theorem \ref{tumb}, there is exists a globally hyperbolic manifold $M_{n}$ diffeomorphic to $S_{n} \times \left [0,1  \right ]$ such that the induced metric on $S_{n} \times \left\{ 1\right\}$ is homotopic to $h_{n}^{+} / \rho^{+}_{n}$ and the induced metric on $S_{n} \times \left\{ 0 \right\}$ is homotopic $h^{-}_{n} / \rho^{+}_{n}$.\\   
The globally hyperbolic manifold $M_{n}$ lifts to a globally hyperbolic convex subset $\Omega_{n}$. Up to normalizing, we can find isometries $V^{\pm}_{n}: (\mathbb{D},h^{\pm}_{n}) \to \partial^{\pm} \Omega_{n}$ such that $\partial V_{n}^{\pm}(p) = (p,p)$ for $p = 0, 1, \infty$. Since there is a unique $\rho_{n}^{+}, \rho_{n}^{-}$ equivariant quasi-symmetric map, it follows that $\Phi_{\Omega_{n}} = f_{n}$.\\
By Proposition \ref{approx-the-quasi}, the maps $f_{n}$ are uniformly quasi-symmetric. Then there is $k' > 1$ such that $\partial_{\infty} \Omega_{n}$ is a $k'$ quasi-circle for any $n$. This implies that $\Omega_{n}$ converge in the Hausdorff topology to a globally hyperbolic convex subset $\Omega$. Then by Proposition \ref{continuty-of-the gluing-maps}, the isometries $V_{n}^{\pm}: (\mathbb{D},h_{n}^{\pm}) \to \partial^{\pm}\Omega_{n}$ converge to an isometry $V^{\pm}: (\mathbb{D},h^{\pm}) \to \partial^{\pm}\Omega$ and the gluing map $\Phi_{\Omega}$ will be equal $f$.
\end{proof}

\bibliographystyle{alpha}
	\bibliography{biblo.bib}
\end{document}